\documentclass[12pt]{article}
\usepackage{amssymb,amsmath,makeidx, amsthm}
\usepackage{enumerate}
\usepackage{psfig}
\usepackage{color}
\usepackage[all]{xy}
\begin{document}
\newtheorem{theorem}{Theorem}[subsection]
\newtheorem{pusto}[theorem]{}
\numberwithin{subsection}{section}
\newtheorem{corollary}[theorem]{Corollary}
\newtheorem{remark}[theorem]{Remark}
\newcommand{\po}{{\perp_{\omega}}}
\newcommand{\mathto}{\mathop{\longrightarrow}\limits}
\def\x{\index}
\def\bp{\begin{pusto}}
\def\ep{\end{pusto}}
\def\bi{\begin{itemize}}
\def\ei{\end{itemize}}
\def\btl{$\blacktriangleleft$}
\def\btr{$\blacktriangleright$}
\def\mn{\medskip\noindent}
\def\n{\noindent}
\def\p{\partial}
\def\eps{\epsilon}
\def\a{\alpha}
\def\d{\delta}
\def\s{\sigma}
\def\ol{\overline}
\def\wt{\widetilde}
\def\Cal{\cal}
\def\Int{\mathrm{Int}\,}
\def\Ker{\mathrm{Ker\,}}
\def\Id{\mathrm{Id}}
\def\id{\mathrm{id}}
\def\su{\mathrm{supp}}
\def\rank{\mathrm{rank}}
\def\dim{\mathrm{dim}}
\def\codim{\mathrm{codim}}
\def\const{\mathrm{const}}
\def\Span{\mathrm{Span}}
\def\span{\mathrm{span}}
\def\inv{\mathrm{inv}}
\def\Hom{\mathrm{Hom}}
\def\Mono{\mathrm{Mono}}
\def\Diff{\mathrm{Diff}}
\def\Emb{\mathrm{Emb}}
\def\emb{\mathrm{emb}}
\def\imm{\mathrm{imm}}
\def\sub{\mathrm{sub}}
\def\Hol{\mathrm{Hol}}
\def\hol{\mathrm{hol}}
\def\Sec{\mathrm{Sec}}
\def\Sol{\mathrm{Sol}}
\def\Supp{\mathrm{Supp}}
\def\bs{\mathrm{bs}}
\def\Vect{\mathrm{Vect}}
\def\Vert{\mathrm{Vert}}
\def\Ver{\mathrm{Ver}}
\def\Clo{\mathrm{Clo}}
\def\clo{\mathrm{clo}}
\def\Exa{\mathrm{Exa}}
\def\Gr{\mathrm{Gr}}
\def\CR{\mathrm{CR}}
\def\Distr{\mathrm{Distr}}
\def\dist{\mathrm{dist}}
\def\Ham{\mathrm{Ham}}
\def\mers{\mathrm{mers}}
\def\Iso{\mathrm{Iso}}
\def\iso{\mathrm{iso}}
\def\tang{\mathrm{tang}}
\def\norm{\mathrm{norm}}

\def\Tang{\mathrm{Tang}}
\def\trans{\mathrm{trans}}
\def\ttr{\mathrm{trans-tang}}
\def\ot{\mathrm{ot}}
\def\loc{\mathrm{loc}}
\def\Symp{\mathrm{Symp}}
\def\symp{\mathrm{symp}}
\def\Lag{\mathrm{Lag}}
\def\isosymp{\mathrm{isosymp}}
\def\subisotr{\mathrm{sub-isotr}}
\def\Cont{\mathrm{Cont}}
\def\cont{\mathrm{cont}}
\def\Leg{\mathrm{Leg}}
\def\Pleat{\mathrm{Pleat}}
\def\isocont{\mathrm{isocont}}
\def\isot{\mathrm{isot}}
\def\isotr{\mathrm{isotr}}
\def\coisot{\mathrm{coisot}}
\def\coreal{\mathrm{coreal}}
\def\real{\mathrm{real}}
\def\comp{\mathrm{comp}}
\def\old{\mathrm{old}}
\def\new{\mathrm{new}}
\def\top{\mathrm{top}}
\def\rG{\mathrm{G}}
\def\CS{\mathrm{CS}}
\def\ND{\mathrm{ND}}
\def\CND{\mathrm{CND}}
\def\fA{\mathfrak {A}}
\def\fa{\mathfrak {a}}
\def\A{\mathcal {A}}
\def\E{\mathcal {E}}
\def\F{\mathcal {F}}
\def\G{\mathcal {G}}
\def\H{\mathcal {H}}
\def\I{\bf {I}}
\def\J{\mathcal {J}}
\def\N{\mathcal {N}}
\def\R{\mathcal {R}}
\def\P{\mathcal {P}}
\def\T{\mathcal {T}}
\def\L{\mathcal {L}}
\def\SS{\mathcal {S}}
\def\U{\mathcal {U}}
\def\D{\mathcal {D}}
\def\V{\mathcal {V}}
\def\W{\mathcal {W}}
\def\Z{\mathcal {Z}}
\def\C{\mathcal {C}}
\def\O{\mathcal {O}}
\def\BH{{\bold H}}
\def\BF{{\bold F}}
\def\BS{{\bold S}}
\def\bbZ{{\mathbb Z}}
\def\bbR{{\mathbb R}}
\def\bbS{{\mathbb S}}
\def\bbC{{\mathbb C}}
\def\bbQ{{\mathbb Q}}
\def\IS{\R_{\sympl-iso}}
\def\IC{\R_{\cont-iso}}
\def\Sp{\mathrm{Sp}}
\def\GL{\mathrm{GL}}
\def\SO{\mathrm{SO}}
\def\Norm{\mathrm{Norm}}
\def\Nor{\mathrm{Nor}}
\def\U{\mathrm{U}}
\def\X{X^{(r)}}
\newcommand{\sign}{\operatorname{sign}}
\def\Op{{\mathcal O}{\it p}\,}
\def\hook{\lrcorner\,}
\def\phi{\varphi}

\newcommand{\corank}{\mathrm{Corank}}
\newcommand{\St}{\mathrm{Stab}}
\newcommand{\wh}{\widehat}
\title{
THE SPACE OF FRAMED FUNCTIONS IS CONTRACTIBLE \\
{\small To Stephen Smale on his 80th birthday}}

\author{
     Y. M. Eliashberg
    \thanks{Partially supported by the NSF grant DMS-0707103}
     \\Stanford University, \\Stanford, CA 94305 USA\\
     {\small eliash@math.stanford.edu}
\and
 N. M. Mishachev
  \\ Lipetsk Technical University, \\
     Lipetsk, 398055 Russia\\
      {\small mishachev@lipetsk.ru}}

\date{}
\maketitle

\begin{abstract}

\n
According to K.~Igusa (\cite{[Ig84]}) a {\it generalized Morse function}
on $M$ is a smooth function $M\to \bbR$ with only Morse and birth-death
singularities and a {\it framed} function on $M$
is a generalized Morse function with an additional structure:  a framing of the 
negative eigenspace at each critical point of $f$. In (\cite{[Ig87]})
Igusa proved that the space of framed generalized Morse functions is $(\dim\,M -
1)$-connected.
J.~Lurie gave  in \cite{[Lu09]}  an algebraic topological proof that the space 
of framed functions is contractible.
In this paper we give a geometric proof of   Igusa-Lurie's theorem using methods 
of our paper \cite{[EM00]}.

\end{abstract}

\tableofcontents

\section{Framed Igusa functions}
\label{s.framed}
\subsection{Main theorem}
\label{ss.main}

This paper is written at a request of D. Kazhdan and  V. Hinich who asked us  
whether we could adjust  our  proof in \cite{[EM00]} of   K. Igusa's  
h-principle for generalized Morse functions  from \cite{[Ig84]} to the case of 
framed generalized Morse functions considered  by   K. Igusa  in his  paper 
\cite{[Ig87]} and more recently by J. Lurie in \cite{[Lu09]}.  We are very happy 
to devote this paper to Stephen Smale whose geometric construction in 
\cite{[Sm58]} plays the central role in our proof (as well as in the proofs of   
many other h-principle type  results.)

\mn Given an $n$-dimensional manifold $W$, a  {\it generalized Morse  function}, 
or as we call it in this paper {\it Igusa 
function},  is  a function 
with only Morse ($A_1$) and birth-death ($A_2$) type singularities. 
A {\it framing} $\,\xi$ of an Igusa  function $\phi:W\to\bbR$ is a 
trivialization of the 
negative eigenspace  of the Hessian quadratic form at $A_1$-points which satisfy 
certain extra conditions at $A_2$-points, see a precise definition
below.

\mn
If  the manifold $W$ is endowed with a foliation $\F$ then 
we call  $\phi:(W,\F)\to\bbR$ a {\it leafwise Igusa}
function if restricted to leaves it has only Morse or birth-death type 
singularities. A {\it framing} $\xi$ of a leafwise Igusa function 
$\phi:(W,\F)\to\bbR$ is a leafwise framing; see   precise definitions
below.

\mn
The following theorem is the main result of the paper. We use Gromov's notation 
$\Op A$ for an unspecified open neighborhood of a closed subset $A\subset W$.
\bp
\label{thm:main}
{\bf (Extension theorem)}
Let $W$ be an $(n+k)$-dimensional manifold with 
an $n$-dimensional foliation $\F$. Let $A\subset W$ be a closed (possibly empty) 
subset and $(\phi_A,\xi_A)$ a framed leafwise Igusa function defined on  $\Op 
A\supset A$. Then there exists  a framed leafwise Igusa function 
$(\phi,\xi)$ on the whole $W$  which coincides with $(\phi_A,\xi_A)$ on   
 $\Op A$.
\ep

Theorem \ref{thm:main} is   equivalent to the fact that the space of 
framed Igusa functions is contractible, which is a content of J. Lurie's 
extension (see Theorem 3.4.7 in \cite{[Lu09]}) of K. Igusa's result  from 
\cite{[Ig87]}. The current form of the theorem allows us  to avoid discussion of  
the topology on this space, comp. \cite{[Ig87]}.

\mn {\it Acknowledgements.} We are grateful to D. Kazhdan and V. Hinich for 
their encouragement to write this paper and to S. Galatius for enlightening 
discussions.

\subsection{Framed Igusa functions}
\label{ss:genuine}
{\it Objects associated with a leafwise Igusa function}.
Let $T\F$ denote the $n$-dimensional subbundle of $TW$ tangent to the 
leaves of the foliation $\F$. Let us fix a Riemannian metric on $W$.
Given  a leafwise Igusa function  (LIF) $\phi$ we associate with it the 
following objects:
\begin{itemize}

\item  $V=V(\phi)$  is    the set of all its leafwise critical points,
i.e. the set of zeros of the leafwise differential $d_\F\phi:W\to T^*\F$. 

\item  $\Sigma=\Sigma(\phi)$   is the set of $A_2$-points.
Generically, $V$ is a $k$-dimensional submanifold of $W$ which is transversal 
to $\F$ at the set  $V\setminus \Sigma$ of $A_1$-points  and has the fold type 
tangency to $\F$ along a $(k-1)$-dimensional submanifold $\Sigma\subset V$ of 
leafwise $A_2$-critical points of $\phi$. 

\item  $\Vert$ is the restriction bundle $T\F|_V$.

\item $d^2_\F\phi$  is the  leafwise quadratic differential of $\phi$. It is 
invariantly defined at  each point $v \in V$. $d^2_\F\phi$ can be  viewed as a 
homomorphism  $\Vert\to\Vert^*$.  Using our choice of a Riemannian metric  we 
identify the  bundles   $\Vert$ and $\Vert^*$  and  view $d^2_\F\phi$  as a 
self-adjoint operator $\Vert\to\Vert$.
This operator  is non-degenerate at the points of $V\setminus \Sigma\,$, and has a 
1-dimensional kernel $\lambda\subset \Vert|_\Sigma$.   Note that  
$\lambda$ is tangent to $V$, and thus we have $\lambda=\Vert\cap TV|_{\Sigma}$. 
 
\item $d^3_\F\phi$ is the invariantly defined  third leafwise differential, which 
is a cubic form on $\lambda$. 
For a leafwise  Igusa  function $\phi$ this cubic form is non-vanishing, and 
hence the bundle $\lambda$ is trivial and can be canonically oriented by 
choosing  the direction in which the cubic function $d^3_\F\phi$ increases. We 
denote by  $\lambda^+$ the unit vector in $\lambda$ which defines its 
orientation. 
 \end{itemize}

\mn 
{\it Decomposition of $\,V(\phi)$ and splitting of $\,\Vert $.} 
The index of the leafwise quadratic differential $d^2_\F\phi\,(v)$, $v\in V$,
may takes values 
$0,1,\dots ,n$ for $v\in V\setminus \Sigma$ and $0,1,\dots ,n-1$ for 
$v\in \Sigma$. Let 
$$V\setminus \Sigma= V^0\cup\dots \cup V^n\,\,\,\,\hbox{and}\,\,\,\,
\Sigma=\Sigma^{1}\cup\dots\cup\Sigma^{n-1}$$
be the decompositions of $V\setminus \Sigma$ and $\Sigma$  
according to the index.
Note that 
$\Sigma^i$ is the intersection of the closures of $V^{i}$ and $V^{i+1}$.
Then for $v\in V^i$ we have 
the splitting 
$$T_v\F=\Vert (v)=\Vert^{i}_+(v)\oplus \Vert^i_-(v)$$
where $\Vert^{i}_+(v)$ and $\Vert^i_-(v)$ are the 
positive and the negative eigenspaces of $d^2_\F\phi(v)$, and for any 
$\sigma\in\Sigma^{i}$ we have the splitting 
$$T_\sigma \F=\Vert(\sigma)=\Ver(\sigma)\oplus 
\lambda(\sigma)\,=
\Ver^i_+(\sigma)\oplus \Ver^i_-(\sigma)\oplus 
\lambda(\sigma)\,$$  
($\Ver\neq\Vert$ !), where $\Ver^i_+(v)$ and $\Ver^i_-(v)$ are the 
positive and the negative eigenspaces of $d^2_\F\phi\,(\sigma)$.
For $\sigma\in \Sigma^i$ and $v\in V^i$ we have 
$$\lim_{v\to\sigma}\Vert^i_+(v)=\Ver^i_+(\sigma)\oplus 
\lambda(\sigma)\,\,\,\,\hbox{and}\,\,\,
\lim_{v\to\sigma}\Vert^i_-(v)=\Ver^i_-(\sigma)\,.$$
For $\sigma\in \Sigma^i$ and $v\in V^{i+1}$ we have 
$$\lim_{v\to\sigma}\Vert^{i+1}_-(v)=\Ver^i_-(\sigma)\oplus 
\lambda(\sigma)\,\,\,\,\hbox{and}\,\,\,
\lim_{v\to\sigma}\Vert^{i+1}_+(v)=\Ver^i_+(\sigma)\,.$$

\begin{figure}[hi]
\centerline{\psfig{figure=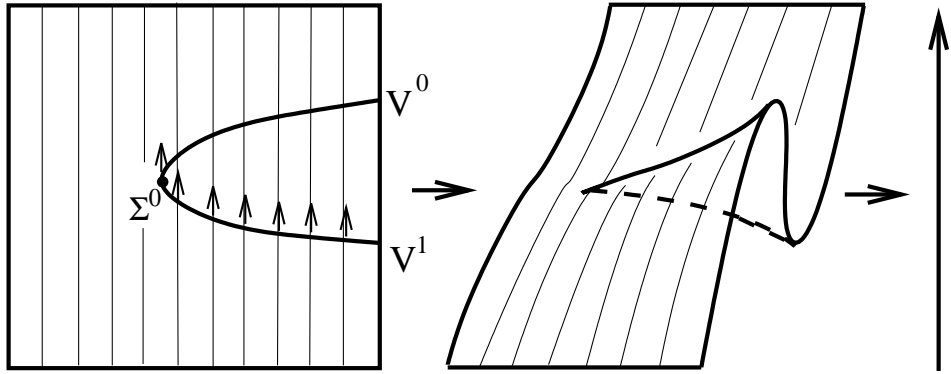,height=50mm}}
\caption {\small Framed leafwise Igusa function}
\label{wi8}
\end{figure}

\mn
{\it Framing of a leafwise Igusa function.}  A framing of a leafwise Igusa 
function $\phi$ is an ordered set 
$\xi=(\xi^1,\dots,\xi^n)$ of unit vector 
fields in  $\Vert(V)$ such that: 
\bi
\item $\xi^i$ is defined (only) over the union 
$\Sigma^{i-1}\cup V^i\cup\dots\cup\Sigma^{n-1}\cup V^n$;
\item $\xi^i|_{\Sigma^{i-1}}=\lambda^+|_{\Sigma^{i-1}}$;
\item $(\xi^1,\dots,\xi^i)|_{V^i}$ is an orthonormal framing for $\Vert_-^i$.
\ei 
In particular, $\,\xi^n$ is defined only on $\Sigma^{n-1}\cup V^n$ and 
$\xi^1$ is defined only on $V\setminus V^0$.
The pair $(\phi,\xi)$ is called a  {\it framed leafwise Igusa function} 
(see Fig.\ref{wi8}).

\mn
The motivation for adding  a framing is discussed in \cite{[Ig87]}.

\subsection{Framed formal leafwise Igusa functions}  
\label{ss:formal}
A {\it formal  leafwise Igusa function}  (FLIF) is a quadruple 
$\Phi=(\Phi^0,\Phi^1, \Phi^2, \lambda^+)$ where: 
\bi
\item $\Phi^0: W \rightarrow \bbR$ is any   function; 
\item $\Phi^1:W\rightarrow T \F\,$ is a  vector field tangent to $\F$,
vanishing on a subset $V=V(\Phi)\subset W$;
\item $\Phi^2$ is a self-adjoint operator   
$\Vert\to\Vert$, which has rank $n-1$
over a subset $\Sigma=\Sigma(\Phi)\subset V$ and rank $n$\, over 
$V\setminus\Sigma$;
\item$\lambda^+$ is a unit vector field  in the line bundle
where   $\lambda:=\Ker(\Phi^2|_{TV|_\Sigma})$. 
\ei
A leafwise $3$-jet of a genuine Igusa function can be viewed as a formal 
Igusa function $\Phi$, where 
$\Phi^0=\phi$, $\Phi^1=\nabla_\F\phi$, $\Phi^2=d^2_\F\phi$ and $\lambda^+ $ is 
the unit vector field in $\Ker d^2_\F$ oriented by the third differential $ 
d^3_\F\phi$. We denote this FLIF $\Phi$ by $J(\phi)$.  A FLIF  $\Phi$ of the 
form $J(\phi)$ is called {\it holonomic}.
Thus we can view a genuine Igusa function as a 
holonomic formal Igusa  function. Usually we will not distinguish between 
leafwise holonomic functions and corresponding holonomic FLIFs.

\mn
Given  a   FLIF $\Phi$ we will use the notation similar to 
  the holonomic case.
Namely,  
  
\bi
\item $V^i\subset V\setminus\Sigma\,$ is the set of points $v\in 
V\setminus\Sigma$ where the index (dimesion of the negative eigenspace) of  
 $\Phi^2_v\,$ is equal to $i\,$, $i=0,\dots, n$;
\item $\Sigma^i\subset\Sigma$ is the set of points $\sigma\in\Sigma$ such that 
the index of $\Phi^2_\sigma\,$ is equal to $i\,$, $i=0,\dots,n-1$;
\item $T_v\F=\Vert (v)=\Vert^{i}_+(v)\oplus \Vert^i_-(v)$ 
where $\Vert^{i}_+(v)$ and $\Vert^i_-(v)$  are the 
positive and the negative eigenspaces of $\Phi^2_v$, $v\in V$;
\item  $T_\sigma \F=\Vert(\sigma)=\Ver(\sigma)\oplus 
\lambda(\sigma)\,=
\Ver^i_+(\sigma)\oplus \Ver^i_-(\sigma)\oplus 
\lambda(\sigma)\,$  where $\Ver^i_+(v)$ and $\Ver^i_-(v)$ are the 
positive and the negative eigenspaces of $\Phi^2_\sigma$,  $\sigma\in 
\Sigma^{i}$.
\ei
As in the holonomic case, for $\sigma\in \Sigma^i$ and $v\in V^i$ we have 
$$\lim_{v\to\sigma}\Vert^i_+(v)=\Ver^i_+(\sigma)\oplus 
\lambda(\sigma)\,\,\,\,\hbox{and}\,\,\,
\lim_{v\to\sigma}\Vert^i_-(v)=\Ver^i_-(\sigma)\,,$$
and so on.

\mn
A {\it framing}  for a formal leafwise Igusa function $\Phi$
is an ordered set $\xi=(\xi^1,\dots,\xi^n)$ of unit vector 
fields in  $\Vert(V)$ such that: 
\bi
\item $\xi^i$ is defined (only) over the union 
$\Sigma^{i-1}\cup V^i\cup\dots\cup\Sigma^{n-1}\cup V^n$;
\item $\xi^i|_{\Sigma^{i-1}}=\lambda^+|_{\Sigma^{i-1}}$;
\item $(\xi^1,\dots,\xi^i)|_{V^i}$ is an orthonormal framing for $\Vert_-^i$.
\ei 
The pair $(\Phi,\xi)$ is called a {\it  framed} formal leafwise Igusa function 
(framed FLIF).  

\begin{figure}[hi]
\centerline{\psfig{figure=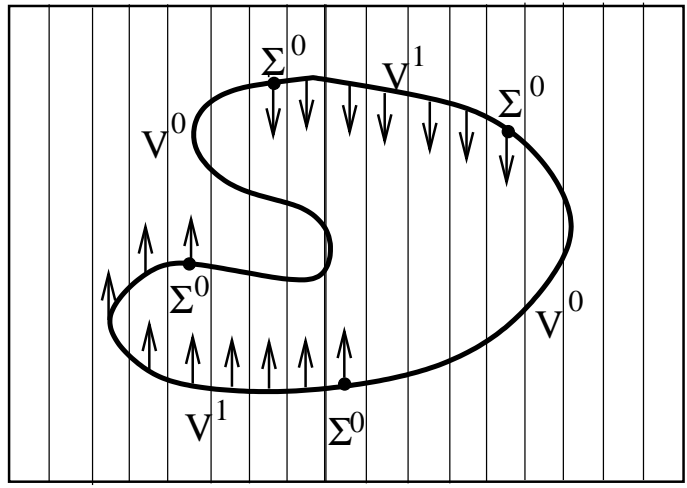,height=65mm}}
\caption {\small Framed FLIF}
\label{wi9}
\end{figure}

\mn
As in the holonomic case, for  a generic FLIF 
$\Phi$ the set  $V$ is a $k$-dimensional manifold and $\Sigma$ its codimension 
1 submanifold. However, $\Sigma$ has nothing to do with    tangency of 
$V$  to $\F$, and moreover there is no control of the type of the tangency 
singularities between $V$ and $\F$ (see Fig.\ref{wi9}).

\mn In what follows we will need to consider FLIFs for different foliations on 
$W$. We will say that $\Phi$ is an $\F$-FLIF when we need to emphasize the 
corresponding foliation $\F$. Moreover, the notion of a FLIF can be generalized 
without any changes to an arbitrary, not necessarily integrable $n$-dimensional 
distribution 
 $\zeta\subset TW$. We will call such an object  a $\zeta$-FLIF.  In the case 
when 
  a distribution $\zeta$  is integrable and integrates into a foliation $\F$ we 
will use as synonyms  both terms: $\zeta$-FLIF  and $\F$-FLIF. 
  
\mn {\it Push-forward operation for FLIFs.}
Let $\zeta ,\wt \zeta $ be two $n$-dimensional distributions in $TW$.
Let $f:W\to W$ be a diffeomorphism covered by an isomorphism $F:\zeta 
\to\wt\zeta\,$.
Let $\Phi $ be  a $\,\zeta $-FLIF. Then we define the push-forward $\wt\zeta $-
FLIF
$\wt\Phi =(f,F)_*\Phi=(\wt\Phi^0,\wt\Phi^1,\wt\Phi^2,\wt\lambda^+)$ as 
  \begin{description}
\item{-} $\wt \Phi ^0(f(x)):=\Phi^0(x),\;x\in W$;
\item{-} $\wt\Phi^1_{f(x)}(F(Z))=\Phi^1_x(Z),\;x\in W,\,Z\in\zeta_x$; 
\item{-} $\wt\Phi^2_{f(x)}(F(Z))=F(\Phi^2_x (Z)),\; x\in V,\,Z\in\Vert_x=\zeta_x\,$;
\item{-}  $\wt\lambda^+(f(x))=F(\lambda^+(x)),\; x\in\Sigma.$
\end{description}
If $\Phi$ is framed   then the push-forward operator $(f,F)_*$ transforms its 
framing $\xi$ to a    framing $\wt\xi$ of   $\wt\Phi$ in a natural way:
 \begin{description}
\item{-} $\wt\xi^i(f(x))=F(\xi^i(x)),\; x\in V$.
\end{description}
\mn Note that if $\zeta $ and $\wt\zeta $ are both integrable, i.e. tangent to 
foliations $\F $ and $\wt\F $,   $F=df$  and  $\Phi$ is holonomic  i.e. 
$\Phi=J(\phi)$     then $\wt\Phi$ is also  holonomic, $\wt\Phi=J(\wt\phi)$,
 where    $\wt\phi =\phi \circ f^{-1}$.

\subsection{Outline of the proof and   plan of the paper}\label{ss.out}

Any framed leafwise Igusa function can be extended from $\Op A$ to $W$ {\it 
formally}, i.e. as a  framed FLIF $(\Phi,\xi)$, see Lemma \ref{thm:formalext2}. 
This is, essentially, an original Igusa's observation from \cite{[Ig87]}.
We then gradually improve $(\Phi,\xi)$ to make it holonomic.
Note that  unlike the holonomic case, the homotopical data   associated 
with $\Phi^1$ and $\Phi^2$ are essentially unrelated. We formulate the 
necessary so-called {\it balancing} homotopical  condition for a FLIF to be
 holonomic, see Section \ref{sec:prep}, and show that one can always 
 make a FLIF  $(\Phi,\xi)$ balanced  via a modification, called {\it stabilization}, see
 Section 
\ref{sec:stab}.

Our next task is to arrange that  $V(\phi)$  has fold type tangency with respect 
to the  
foliation $\F$, as it is supposed to be in the holonomic case. FLIFs satisfying 
this property, together with certain  additional coorientation conditions over 
the fold, are 
called {\it prepared}, see Section \ref{ss:notation}. We   observe that for  a 
prepared FLIF one can define a stronger    necessary  homotopical condition for 
holonomicity. We call prepared  FLIFs satisfying this stronger condition {\it 
well balanced}, see Section \ref{sec:prep}.

Given any FLIF $\Phi$ one can associate with it a {\it twisted   normal bundle} 
(also called   {\it virtual vertical bundle}) ${}^\Phi\Vert$ over $V=V(\Phi)$ which is a subbundle of $TW|_{V}$ obtained by 
twisting the normal bundle  of $V$ in $W$ near $\Sigma=\Sigma(\Phi)$, \,see 
Section \ref{sec:twisted}. In the holonomic case we have ${}^\Phi\Vert=\Vert$,
 see \ref{lm:hol-prep}. A crucial observation is that  the manifold $V$ has fold type 
tangency to any extension $\zeta$ of the bundle ${}^\Phi\Vert$ to a neighborhood 
of $V$, see \ref{lm:tang-to-twist}. Moreover, if $\Phi$ is balanced then there 
exist a global extension $\zeta$ of ${}^\Phi\Vert$  and a bundle isomorphism 
$F:\Vert\to{}^\Phi\Vert$ homotopic  to the  identity $\Vert\to\Vert$ through 
injective bundle homomorphisms into $TW|_V$ such that the push-forward framed  
$\zeta$-FLIF $(\wt\Phi,\wt\xi)=(\Id,F)_*(\Phi,\xi)$ is well balanced, see 
\ref{lm:norm-prep}. 
If $\Phi$ is holonomic on $\Op A$ then the bundle $\zeta$ and the framed  FLIF 
$(\wt\Phi,\wt\xi)$ coincide with $T\F$ and $(\Phi,\xi)$ over $\Op A$.  

The  homotopy of the homomorphism $F$ generates a homotopy 
of distributions $\zeta_s$ connecting $\zeta$ and $T\F$. 
If it were possible to  construct a fixed on $\Op A$ isotopy $V_s$ of $V$ in $W$ 
keeping $V_s$ folded with respect to  $\zeta_s$ then one could cover this 
homotopy by a  fixed on $\Op A$ homotopy of  framed well balanced $\zeta_s$-
FLIFs $(\wt\Phi_s,\wt\xi_s)$ beginning with  
$(\wt\Phi_0,\wt\xi_0)=(\wt\Phi,\wt\xi)$.  Though this is, in general, 
impossible, the wrinkling embedding theorem from \cite{[EM09]} allows us to do 
that after a certain additional modification of $V$, called {\it pleating}, see 
Theorem \ref{thm:folds}. We then show that the pleating construction can be 
extended 
to the class of   framed well balanced  FLIFs, see 
Section \ref{sec:pleat-FLIF}. Thus we get a framed well balanced  FLIF $(\wh 
\Phi,\wh\xi)$ extending  the local  framed leafwise Igusa function 
$(\phi_A,\xi_A)$.

The proof now is concluded in two steps. First, we show, see Lemma
\ref{lm:loc-int}, that a framed  well balanced FLIF can be made holonomic near 
$V$, and then use the wrinkling theorem from \cite{[EM97]} to construct  a 
holonomic extension to the whole manifold $W$,   see Step 5 in Section 
\ref{sec:proof}.

\mn  The paper has the following organization.
  In Section \ref{sec:distr-fold}  
 we discuss the notion of fold tangency of  a submanifold with respect to a not 
necessarily integrable distribution, define the pleating construction for 
submanifolds and formulate the main  technical  result, 
Theorem  \ref{thm:folds}, which
  is an analog  for folded maps of Gromov's directed embedding theorem, see 
\cite{[Gr86]}. This is a 
corollary of the results of  \cite{[EM09]}.  
  Section \ref{sec:FLIF} is the main part of the paper. We define and study 
there the 
notions  and properties of balanced, prepared and well balanced FLIFs, and 
gradually realize the 
described above program of making a framed FLIF well balanced, see Proposition 
\ref{sec:FLIF}.   We also prove here Igusa's result about existence of a formal 
extension for framed FLIFs, see \ref{thm:formalext2}, and   local 
integrability of well balanced FLIFs, see \ref{lm:loc-int}. Finally, in Section 
\ref{sec:proof} we just recap the main steps of the proof.

 \section{Tangency of a submanifold to a distribution}\label{sec:tangency}
 In   this section we always denote by   $V $   an $n$-dimensional submanifold 
of an $(n+k)$-dimensional manifold $W$,  by $\Sigma$ a  codimension 1 
submanifold of $V$ and  by $\Norm=\Norm(V)$ the normal bundle of $V$.
 
 \subsection{Submanifolds folded with respect to a 
distribution}\label{sec:distr-fold}
 Let $\zeta$ be  an $n$-dimensional distribution, i.e. a subbundle 
$\zeta\subset TW$.
  The {\it non-transversality} condition of $V$ to $\zeta$ defines a variety 
$\Sigma_\zeta$ of the $1$-jet space $J^1(V,W)$. 
We say that $V$ has  at a point $p\in V$  a tangency to $\zeta$ of 
{\it fold} type  if 
\begin{itemize}
\item $\corank\, \pi^{\,\zeta}|_{T_pV}=1$;
\item $J^1(j):V\to J^1(V,W)$, where $j:V\hookrightarrow W$ is the inclusion, is 
transverse to $\Sigma_\zeta$; we denote $\Sigma:=(J^1(j))^{-1}(\Sigma_\zeta)$;
\item $\pi^{\,\zeta}|_{T_p\Sigma }:T_p\Sigma \to TW_{p}/\zeta$
is injective.
\end{itemize}

If $\zeta$ is integrable, and hence locally is tangent to an affine foliation 
defined by the projection $\pi:\bbR^{n+k}\to \bbR^k$, these conditions are  
equivalent to the 
requirement that the restriction $\pi|_V$ has fold type singularity, and in this 
case one has a normal form for the fold tangency.

\mn If $V$ has fold type tangency to $\zeta$ along $\Sigma$ then we say that $V$ 
is 
{\it folded} with respect to $V$ along $\Sigma$.
The fold locus $\Sigma\subset V$ is a codimension one submanifold, and at each 
point $\sigma\in\Sigma$ the  $1$-dimensional line field $\lambda=\Ker 
\pi^\zeta|_{TV}=\zeta|_{V}\cap\, TV$ is transverse to $\Sigma$.

\begin{figure}[hi]
\centerline{\psfig{figure=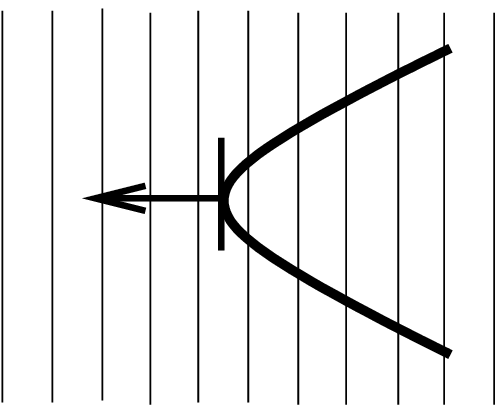,height=30mm}}
\caption {\small Characteristic coorientation of the fold}
\label{nwi31}
\end{figure}

\mn The hyperplane field $T\Sigma\oplus\zeta|_{\Sigma}$ can be canonically
 cooriented. 
 In the case when $\zeta$ is integrable    this  coorientation can be defined 
as follows.
  The leaves of the foliation trough points of $ \Sigma$ form a hypersurface 
which divides a sufficiently small tubular  neighborhood $\Omega$ of $ \Sigma$ 
in $W$ into two parts, $\Omega=\Omega_+\cup\Omega_-$, where    $\Omega_-$ is the 
part which contains $V\cap\Omega$.
Then   the  {\it characteristic coorientation}  of the fold  $\Sigma$ is the   
coorientation of   the hyperplane $T\Sigma\oplus\zeta|_{\Sigma}$      
determined by  the  outward normal vector field to $\Omega_-$ along $ \Sigma $,
see Fig.\ref{nwi31}.
 For  a general $\zeta$     take a point $\sigma_0\in\Sigma$, a neighborhood $U$ 
of $\sigma_0$ in $V$, an  arbitrary unit vector  field 
$\nu^+\in (T\Sigma\oplus\zeta|_{\Sigma})^\bot$ and  
consider an embedding 
$g:U\times(-\eps,\eps)\to W$ such that $g(x,0)= x,\, x\in U$ and 
$\frac{\p g}{\p t}(\sigma, 0)=\nu^+$, $\sigma\in\Sigma\cap U$,
where $t\in(-\eps,\eps)$  is the coordinate corresponding to the second factor.
Consider the line field $L=d\,g(\zeta)$. Note  that 
 $L|_{(\Sigma\cap U)\times 0}=\lambda$.
The line field $L$ integrates to a $1$-dimensional foliation on $U\times (-
\eps,\eps)$ which has   a  tangency of fold type to $U\times 0$ along 
$\Sigma\times 0$. Hence, $U\times 0\subset U\times(-\eps,\eps)$ can be
 cooriented, as in the integrable case, which gives the required coorientation
 of $ T\Sigma\oplus\zeta|_{ \Sigma }$, see Fig.~\ref{nwi31}.
 
 \mn It is important to note that {\it the property that $V$ has a  fold type 
tangency   to $\zeta$ along $\Sigma$ depends only on $\zeta|_V$, and not on its 
extension to $\Op V$.}
Similarly,  the above definition of the characteristic coorientation  of 
$ T\Sigma\oplus\zeta|_{ \Sigma }$ is independent of all the choices and depends 
only on $\zeta|_V$ and not on its extension to $\Op V$.

\mn The following simple lemma (which we do not use in the sequel) clarifies the 
geometric meaning of the fold tangency.

\bp\label{lm:fold-normal}{\bf (Local normal form for fold type tangency to a
distribution)} 
Suppose $V\subset W$ is folded with respect to $\zeta$ 
along $V$ and the fold  $\Sigma$ is cooriented. Denote  
$\lambda:=\zeta|_\Sigma\cap TV|_\Sigma$ and $\eta:=(\zeta|_V)/\lambda$. Consider 
the pull-back  $\wt\eta$ of the bundle $\eta$ to $\Sigma\times\bbR^2$ and denote 
by $E$ the total space of this bundle.   Then there exists a neighborhood   
$\Omega$ of $\Sigma\times 0$ in  $E$, a neighborhood 
$\Omega'\supset\Sigma$ in $W$,   and a diffeomorphism $ \Omega\to\Omega'$  
introducing coordinates $(\sigma,x,z,y)$ in $\Omega'$, $\sigma\in\Sigma\,,\, 
(x,z)\in\bbR^2$, $y\in\eta$, such that 
in these coordinates the manifold $V$ is given by the equations
$z=x^2,y=0$ and the bundle $\zeta|_V$ coincides with the restriction to $V$ of    
the projection $(\sigma,x,z,y)\to(\sigma,z)$. 
 \ep

 Lemma \ref{lm:fold-normal} implies, in particular, that if $V$ is folded with 
respect to $\zeta$ then  $\zeta|_V$ always 
admits an integrable extension to a neighborhood of $V$.

\subsection{Pleating}\label{sec:pleating}
Suppose $V$ is folded with respect to $\zeta$ along $\Sigma$. Let
  $S\subset V\setminus\Sigma$ be a closed codimension $1$ 
   submanifold and      $\nu^+\in\zeta$  be a vector field  defined over $\Op 
S\subset W$. For a sufficiently small $\eps,\delta>0$ there exists an embedding  
$g:S\times[-\delta,\delta]\times[-\eps,\eps]\to W$ such that 
   \begin{itemize}
   \item $\frac{\p g}{\p u}(s,t,u)=
\nu^+(g(s,t,u)),\; (s,t,u)\in S\times[-\delta,\delta]\times [-\eps,\eps]$,
\item 
$g|_{S\times 0\times 0}$ is the inclusion $S\hookrightarrow V$,
   \item $g|_{S\times[-\delta,\delta]\times0}$ is a diffeomorphism onto the 
tubular $\delta$-neighborhood $U\supset S$ in $V$, which sends intervals  
$s\times[-\delta,\delta]\times 0$, $s\in S$, to geodesics normal to $S$.
\end{itemize}

\begin{figure}[hi]
\centerline{\psfig{figure=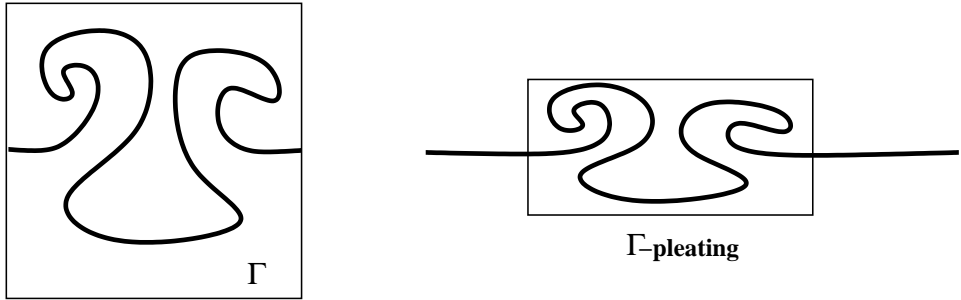,height=40mm}}
\caption {\small $\Gamma$-pleating}
\label{nwi3}
\end{figure}

Let $\Gamma\subset  P:=[-1,1]\times[-1,1]$ be an embedded connected curve which 
near 
$\p P$ coincides with  the line $\{u=0\}$. Here we denote by $t,u$ the 
coordinates corresponding to the two factors.
 We assume that $\Gamma$ is folded with respect to the foliation defined by the 
projection $(t,u)\mapsto t$ (this is a generic condition). We denote by $ 
\Gamma_{\delta,\,\eps}$ the image of 
$\Gamma$ under the scaling $(t,u)\mapsto (\delta t,\eps u)$.
  Consider  a manifold $\wt V$  obtained from $V$ by replacing the   
neighborhood $U$ by a deformed  neighborhood 
  $\wt U_\Gamma=g(S^{n-1}\times\Gamma_{\delta,\,\eps})$.
   We say $\wh V$ is  the result  of   {\it $\Gamma$-pleating}  of $V$ over $S$ 
in the direction of the vector field $\nu^+$, see Fig. \ref{nwi3}.
   
\begin{figure}[hi]
\centerline{\psfig{figure=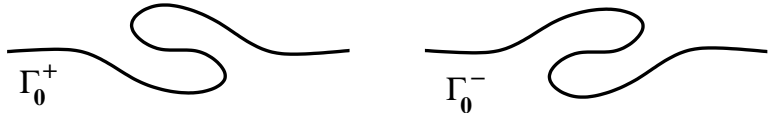,height=15mm}}
\caption {\small  Curves $\Gamma_0^\pm$}
\label{nwi4}
\end{figure}

 The  $\Gamma^+_0$-pleating  with the curve $\Gamma^+_0$ shown on  Fig. 
\ref{nwi4} 
will be referred  simply as   {\it pleating}.  

 \bp\label{thm:folds}
 {\bf (Pleated isotopy)}
 Suppose $V\subset (W,\zeta)$ is folded with respect  to $\zeta$ along 
$\Sigma\subset V$. Let $\zeta_s$, $s\in[0,1]$, be a family of $n$-dimensional 
distributions over a neighborhood   $\Omega\supset V$. Then there exist
\begin{description}

\item{-}  a manifold $\,\wt V\subset \Omega$ obtained from $V$ by a sequence of 
pleatings over boundaries of small  embedded balls in the direction of vector 
fields which extend to these balls, and  

\item{-}  a $\,C^0$-small isotopy   $h_s:\wt V\to \Omega$, 
\end{description}
such that for each $s\in[0,1]$ the manifold  $h_s(\wt V)$ has only  fold type 
tangency to $\zeta_s$. If $\wt\Sigma=\Sigma\cup \Sigma\,'$ is the fold of $\wt V$ 
with respect to $\zeta_0$ then $h_s(\wt\Sigma)$ is the fold    of $h_s(\wh V)$  
with respect to  $\zeta_s$. If the  homotopy $\zeta_s$ is fixed over a neighborhood $Op A$ of
a closed subset $A\subset V$ then one can arrange that $V\cap \Op A=\wt V\cap\Op 
A$ and  the isotopy $h_s$ is fixed over $\Op A$.     \ep

Theorem \ref{thm:folds} is a version of the wrinkled embedding theorem from 
\cite {[EM09]}, see Theorem 3.2 in \cite {[EM09]}
 and the discussion  in  Sections 3.2 and   3.3  in that paper on how to 
  replace the wrinkles 
by spherical double folds  and  how to  generalize Theorem 3.2 to the case of 
not necessarily integrable distributions.
Another cosmetic difference  between the formulations in \cite {[EM09]} and 
Theorem \ref{thm:folds} is that the former one allows not only double folds, but 
also their embryos, i.e. the moments of death-birth of double folds. This can be 
remedied by preserving the double folds till the end in the near-embryo state, 
rather than killing them, and similarly by creating the necessary    number of 
folds by pleating at  the necessary places before the deformation begins.

 \bp\label{rem:special-Gamma}{\bf (Remark)} {\rm
If  $\wt V$ satisfies the conclusion of Theorem \ref{thm:folds} then any  
manifold  $\wt{\wt V}$ obtained from $\wt V$ by an  additional $\Gamma$-pleating 
with {\it any} $\Gamma$ will also have this property.
For our purposes we will  need to     pleat  with  three special curves
 $\Gamma_1$ and $\Gamma_2^\pm$   shown on Figure \ref{nwi5}.
 \begin{figure}[hi]
\centerline{\psfig{figure=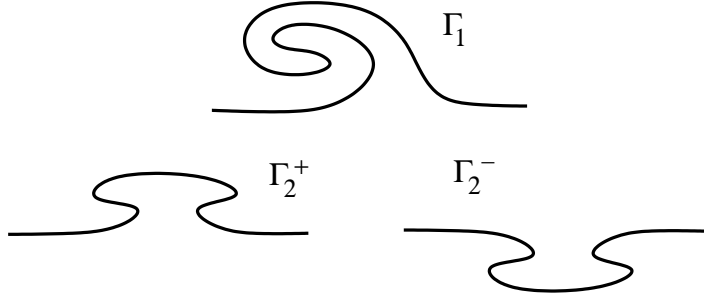,height=40mm}}
\caption {\small  Curves $\Gamma_1$ and $\Gamma^\pm_2$}
\label{nwi5}
\end{figure}
As it clear from this picture,   a   pleating with any of these curves  can be 
viewed as a result of    
a $\Gamma_0^+$-pleating followed by a second  $\Gamma_0^-$-pleating. Hence, in 
the formulation of Theorem 
\ref{thm:folds} one can  pleat with any of the curves  $\Gamma_1$ and 
$\Gamma^\pm_2$  instead of $\Gamma^+_0$.}
\ep

\section{Geometry of FLIFs}\label{sec:FLIF}

\subsection{Homomorphisms $\Gamma_\Phi$ and $\Pi_\Phi$}
\label{ss:notation}

\mn  Given  a  $\zeta$-FLIF $ \Phi$ we will associate with it several objects 
and constructions.

\mn {\it Isomorphism $\Gamma_\Phi: \Norm\to \Vert$.} This isomorphism is 
determined by $\Phi^1$. The tangent bundle $T(\zeta|_V)$ to the the total space 
of the bundle $\zeta|_V$ 
  canonically splits as $\Vert\oplus TW|_V$, 
and hence the  bundle of tangent planes to the section $\Phi^1$ along its 
$0$-set $V$ can be viewed as  a graph  of a  homomorphism $ 
\wh\Gamma_\Phi:TW|_V\to\Vert$  vanishing on $TV$.
 The restriction of this homomorphism to $\Norm$ will be denoted by 
$\Gamma_\Phi$.
The transversality of the section $\Phi^1$ to the $0$-section ensures that  
$\Ker\,\wh\Gamma_\Phi=TV$ and hence $\Gamma_\Phi$ is an isomorphism. 

\begin{figure}[hi]
\centerline{\psfig{figure=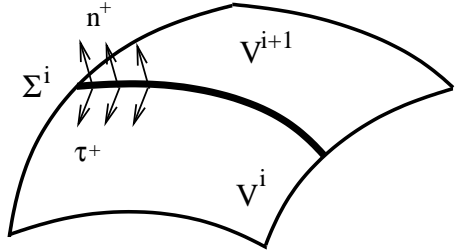,height=35mm}}
\caption {\small The vector fields $n^+$ and $\tau^+$}
\label{nwi6}
\end{figure}

\mn
By an {\it index} coorientation of $\Sigma$ in $V$  we will mean its 
coorientation by 
a normal vector field  $\tau^+$ pointing in the direction of   {\bf decreasing}  
of the index, i.e. on $\Sigma^{\,i}$ it points into $V^i$.
We will denote  by $n^+$ the vector field $\Gamma_\Phi^{-
1}(\lambda^+)\in \Norm(V)\,$, see Fig. \ref{nwi6}.

In the holonomic situation the  index coorientation is given by the vector 
field 
$\lambda^+$ and the vector field $n^+$ determines the characteristic 
coorientation of the fold.

\mn We call a  $\zeta$-FLIF $ \Phi  $ {\it prepared} if 
 \begin{itemize}
 \item $V(\Phi)$ is folded with respect to $\zeta$ with the fold along  
$\Sigma(\Phi)$;
\item  $TV\cap\Vert_{\Sigma}=\lambda$ and  the vector field $\lambda^+$ 
determines the  {\it index} coorientation of the fold;
 \item  the vector field $n^+ $ determines the {\it characteristic} coorientation of 
the fold $\Sigma$.
 \end{itemize}

  Thus, any {\it holonomic}  FLIF (when, in particular, $\zeta$ is integrable)
is prepared.

\mn
{\it Isomorphism $\Pi_\Phi:\Norm\to\Vert$}.
 Given a {\it prepared} $\zeta$-FLIF $\Phi$, let us denote by $K $ the 
restriction of the orthogonal projection $TW|_V\to\Norm$ to the subbundle 
$\Vert=\zeta|_V\subset TW|_V$. The homomorphism $K$ is non-degenerate over 
$V\setminus  \Sigma$ and has  a $1$-dimensional kernel $\lambda$ over $\Sigma$.

 \bp\label{lm:Delta} 
 {\bf (Definition of $\Pi_\Phi$)}
 The composition $\Phi^2\circ K^{-
1}:\Norm|_{V\setminus\Sigma}\to\Vert|_{V\setminus\Sigma}$ continuously extends 
to a non-degenerate homomorphism
 $\Pi_\Phi:\Norm\to\Vert$.
 \ep 
 {\it Proof.}  
Let us  prove the extendability of the inverse operator 
$K\circ\left(\Phi^2\right)^{-1}$.
There exists a canonical extension of the vector   field $\lambda^+$
as a unit  $\Phi^2$-eigenvector field  $\wt\lambda^+$ on $\Op\Sigma\subset V$. 
 Then 
$$\Phi^2 (\wt\lambda^+(v))=c(v)\wt\lambda^+(v)\,,\; v\in\Op\Sigma\,,$$
where the 
eigenvalue function $c:\Op\Sigma\to\bbR$ has $\Sigma$ as its regular $0$-level. 
Denote
 $\wt\Ver=\wt\lambda^\perp(v)$ the orthogonal eigenspace of $\Phi^2(v)$. Denote 
$\wt\Nor:=K(\wt\Ver)$. The  operator 
$K\circ\left(\Phi^2\right)^{-1}$ is well defined on  
$\Ver=\wt\Ver|_\Sigma\subset\Vert|_\Sigma$ and maps it isomorphically onto 
$\Nor=\wt\Nor|_\Sigma$.  It remains  to prove existence of a non-zero limit
 $$\lim\limits_{v\to v_0\in\Sigma}K\left(\left(\Phi^2\right)^{-
1}(\wt\lambda^+)\right)=\lim\limits_{v\to 
v_0\in\Sigma}\frac1{c(v)}K(\wt\lambda^+(v))\,.$$
The vector-valued function $K(\wt\lambda^+(v))$ vanishes on $\Sigma$ while
the function $c(v)$ has no critical points on $\Sigma$. Hence, the above limit 
exists. On the other hand,  the transversality condition for the fold implies 
that $||K(\wt\lambda^+(v))||\geq a\,\dist(v,\Sigma)$, while $|c(v)|\leq 
b\,\dist(v,\Sigma)$  for some positive constants $a,b>0$, and therefore 
$\lim\limits_{v\to v_0\in\Sigma}\frac1{c(v)}K(\wt\lambda^+(v))\neq 0.$
\qed

\bp\label{lm:hol-bal2}
If  $\Phi$ is holonomic then $\Pi_\Phi=\Gamma_\Phi$.
\ep
{\bf Proof.}  Indeed, recall that $\Gamma_\Phi=\wh\Gamma_\Phi|_{\Norm}$, 
where $\wh\Gamma_\Phi:TW|_V\to\Vert$ is the homomorphism defined by the section 
$\Phi^1$ linearized along its zero-set $V$.
In the holonomic situation   one has over\, $ V\setminus\Sigma$ the equality
$$\wh\Gamma_\Phi|_\Vert=d^2\phi=\Phi^2\,,$$
where $\phi=\Phi^0$. But $\wh\Gamma_\Phi|_\Vert$ and 
     $\wh\Gamma_\Phi|_\Norm$ are related by a projection along the kernel of 
$\wh\Gamma_\Phi$ which is equal to $TV$. Hence,    $\Gamma_\Phi= \Phi^2\circ 
K^{-1}=\Pi_\Phi$. By continuity, the equality    $\Pi_\Phi=\Gamma_\Phi$ 
holds everywhere.      \qed

 \subsection{Twisted normal bundle and the isomorphism 
$\Delta_\Phi$}\label{sec:twisted} 
Given any $\zeta$-FLIF $\Phi$ we define  here a {\it twisted normal bundle},
 or as we also call it {\it virtual vertical bundle} ${}^\Phi\Vert \subset 
TW_V$ over $V$. As we will see later (see \ref{lm:hol-prep}),
 in the holonomic  case  ${}^\Phi\Vert $   coincides with $\Vert$.  

\mn 
Let  $U=\Sigma\times[-\eps,\eps]$ be the tubular neighborhood of  $\Sigma$ in 
$V$  of radius $\eps>0$. We assume that the splitting is chosen in such a way 
that the vector field $\frac{\p}{\p t}$, where $t$ is the coordinate 
corresponding to the second factor, defines the index coorientation of $\Sigma$ 
in $V$, and hence coincides with $\tau^+$. We denote  
$$U_+:=\Sigma\times(0,\eps],\,\,\,U_-:=\Sigma\times[-\eps,0)\,.$$
Denote by 
$\wt\lambda^+\in\Vert$ the unit eigenvector field of $\Phi^2|_U$ which extends 
$\lambda^+\in\Vert|_\Sigma$.
 If $\eps$ is small enough then such extension is uniquely defined. Let 
$\wt\Ver:=\wt\lambda^\perp$ be the complementary $\Phi^2$-eigenspace. We have 
$\wt\Ver|_\Sigma=\Ver$.
Denote $\wt n^+:=\Gamma_\Phi^{-1}(\wt\lambda^+)$,\,\,
$\wt\Nor=\Gamma_\Phi^{-1}(\wt\Ver),\,\, \wt\tau^+:=\frac{\p}{\p t}$.
Choose a function $\theta:U\to [-\frac\pi2,\frac\pi2]$ which  has $\Sigma$ as 
its regular level set   $\{\theta=0\}$,  and which is  equal to 
$\pm\frac\pi2$  near $\Sigma\times(\pm\eps)$.   
  
\begin{figure}[hi]
\centerline{\psfig{figure=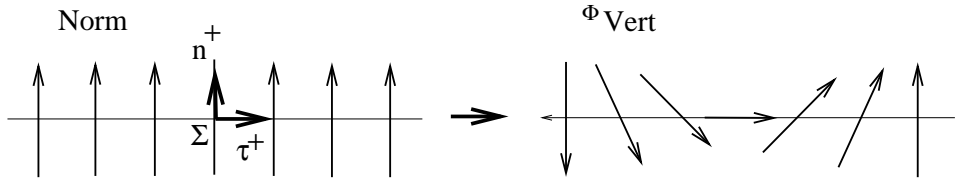,height=25mm}}
\caption {\small Twisting   the normal bundle}
\label{nwi7}
\end{figure}

\mn
We define the bundle ${}^\Phi\Vert$ in the following way.  Over $V\setminus 
U$ it is    equal to $\Norm$. The fiber over a point  $v\in U$  is equal to
   $\Span(\wt\Nor,\mu(v))$, where the line
 $\mu(v)$ is generated by the vector
$$\mu^+(v)=\sin    \theta(v)  \wt n^+(v)+
 \cos   \theta(v)\wt \tau^+(v),$$
 see Fig. \ref{nwi7}.

\mn
{\it Isomorphism $\Delta_\Phi:\Vert\to{}^\Phi\Vert$}.
Let $c:U\to\bbR$ be the eigenvalue function corresponding to the 
$\Phi^2$-eigenvector field $\wt\lambda^+\in\Vert$ on $U$, i.e. we have 
$\Phi^2(\wt\lambda^+(v))=c(v)\wt\lambda^+(v)$, $v\in U$. The function $c$  is 
positive on $U_+$ and negative on $U_-$.
Let $\wt c:U\to\bbR$ be  any {\it positive } function which is equal to $c$ on 
$\p U_+=\Sigma\times\eps$ and equal to $-c$ on  $\p U_-=\Sigma\times(-\eps)$.
 
  We then   define the operator $$\Delta_\Phi:\Vert\to
 {}^\Phi\Vert$$ by the formula

\begin{equation}\label{eq:Delta}
\Delta_\Phi(Z)= \begin{cases}
   \Gamma_{\Phi}^{-1}(\Phi^2(Z)),&\hbox{over}\; V\setminus U, \; Z\in\Vert\,;  
\cr
    \Gamma_{\Phi}^{-1}(\Phi^2(Z)),&\hbox{over}\;  
\;U,\;  Z\in\wt\Ver\,;  \cr
   \wt c(v)\left(\sin    \theta(v)  \wt n^++
 \cos   \theta(v) \wt\tau^+\right) ,&Z=\wt\lambda^+(v),\;
  v\in U\,.
   \end{cases}
\end{equation}

 \mn  

It will be convenient for us to keep some ambiguity  in the definition of 
${}^\Phi\Vert$ and $\Delta_\Phi$. However, we note that the space of choices we 
made in the definition is contractible, and hence the objects are defined in a 
homotopically canonical way.

\mn
Let us extend ${}^\Phi\Vert$ and   $\Delta_\Phi$   to a 
neighborhood $\Op V\subset W$. We will keep the same notation for the extended 
objects.

\bp\label{lm:tang-to-twist}
 For any $\zeta$-FLIF $\Phi$  the  ${}^\Phi\Vert$-FLIF 
$$\Phi^{\Norm}=(\Id,\Delta_\Phi)_*\Phi$$
on $\Op V$ is prepared.
\ep

\begin{figure}[hi]
\centerline{\psfig{figure=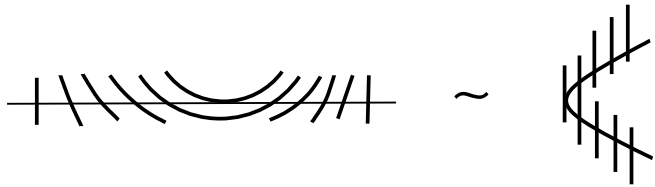,height=25mm}}
\caption {\small $V$ is folded with respect to ${}^\Phi\Vert$.}
\label{nwi9}
\end{figure}

{\bf Proof.} Denote $\wh\Phi:=\Phi^\Norm$.  We have $ V(\wh\Phi)=V(\Phi)=V$.
First of all we observe (see Fig.~\ref{nwi9}) that 
$V$ is folded with respect to ${}^\Phi\Vert$
  along $\Sigma$ and the vector field  
$n^+=n^+(\Phi)$ defines  the characteristic coorientation of the fold. On the 
other hand, 
$\lambda^+(\wh\Phi)=\Delta_\Phi(\lambda_+(\Phi))=\tau^+(\Phi)=\tau^+(\wh\Phi)$ 
and  
 $$n^+(\wh\Phi)=\Gamma_{\wh\Phi}^{-1}(\lambda^+(\wh\Phi))=\Gamma^{-
1}_{\wh\Phi}(\tau^+)=\Gamma_\Phi^{-1}(\Delta_\Phi^{-1}(\tau_+))=\Gamma_\Phi^{-
1}(\lambda^+)=n^+(\Phi),$$ and hence $n^+(\wh\Phi)$ defines the 
characteristic coorientation of the fold $\Sigma$. Thus $\wh\Phi$ is prepared.
\qed

\bp\label{lm:2-bal} For any FLIF $\Phi$ the diagram
$$
\xymatrix{\Norm \ar[rd]_{\Pi_{\wh \Phi}}\ar[rr]^{\Gamma_\Phi}& 
&\Vert\ar[ld]^{\Delta_\Phi} \\&{}^\Phi\Vert&}
 $$
  commutes for appropriate choices in the definition of $\,{}^\Phi\Vert$ and 
$\,\Gamma_\Phi$.
\ep
{\bf Proof.} We need to check that $\Pi_{\wh\Phi}=\Delta_\Phi\circ\Gamma_\Phi$. 
First, we check the equality over $V\setminus U$. We have 
$\Norm|_{V\setminus U}={}^\Phi\Vert|_{V\setminus U}$, and hence $K_{\wh\Phi}=\Id.$ 
Furthermore, over $V\setminus U$ we have 
$\Pi_{\wh\Phi}=K_{\wh\Phi}^{-
1}\circ\wh\Phi^2=\Delta_\Phi\circ\Phi^2\circ\Delta_\Phi^{-1}=
\Gamma_\Phi^{-1}\circ\Phi^2\circ\Phi^2\circ\left(\Phi^2\right)^{-1}\Gamma_\Phi=
\Gamma_\Phi^{-1}\circ\Phi^2\circ\Gamma_\Phi=\Delta_\Phi\circ\Gamma_\Phi.$  
Similarly, we check that 
$\Pi_{\wh\Phi}|_{\wt\Nor}=\Delta_\Phi\circ\Gamma_\Phi|_{\wt\Nor}$.
Finally,   evaluating both parts of the equality on the vector field $n_+$ we 
get:
$\Pi_{\wh\Phi}(n^+)=\lambda^+=\Delta_\Phi(\Gamma_\Phi(n^+))$. Then this implies 
$\Pi_{\wh\Phi}(\wt n^+)=\Delta_\Phi(\Gamma_\Phi(\wt n^+))$ for an appropriate 
choice of the function $\wt c>0$ in the definition of the homomorphism 
$\Delta_\Phi$.
\qed

\subsection{The holonomic case}\label{sec:normal-form}

  We will need the following normal form for   a  leaf-wise Igusa function 
$\phi$ near $\Sigma$  (see \cite{[Ar76], [El72]}).

\mn
  Consider the pull-back of the bundle  $\Ver=\Ver_+\oplus\Ver_-$ defined over 
$\Sigma$ to   $\Sigma\times\bbR\times\bbR$  via  the projection    
$\Sigma\times\bbR\times\bbR\to\Sigma$.  Let $E$ be the total space of this 
bundle. The submanifold $\Sigma\times0\times0 $ of the $0$-section of this 
bundle we will denote simply by $\Sigma$. Consider a function
  $\theta:E\to\bbR$ given by the formula
\begin{equation}\label{eq:cusp}\theta(\sigma,x,z,y_+,y_-)=x^3-
3zx+\frac12(||y_+^2||-||y_-^2||)\,; \end{equation}
  $(\sigma,x,z)\in \Sigma\times\bbR\times\bbR,\;  y_\pm \in(\Ver_\pm)_{\sigma}. 
$
     
Consider the projection  $p:E\to\Sigma\times\bbR$  defined by the formula
$$p(\sigma,x,z,y_+,y_-)=(\sigma,z).$$
There  exists an embedding $g:\Op\Sigma\to W$, where $\Op\Sigma$ is a 
neighborhood of $\Sigma$ in $E$, such that
 \begin{itemize}
 \item $g(\sigma) =\sigma, \sigma\in\Sigma$;
 \item $g$ maps the fibers of the projection $p$ to the leaves of the foliation 
$\F$.
   \item $\phi\circ g=\theta$.
    \end{itemize} 
    Via the parameterization map $g$ we will view $(\sigma, x,z,y_+,y_-)$ as 
coordinates
in $\Op\Sigma\subset W$. In these coordinates the function $\phi$ has the form 
\eqref{eq:cusp}, the manifold $V$ is given by the equations $z=x^2,y_\pm=0$, the 
foliation $\F$ is given  by the fibers of the projection $p$,
the vector field $-\frac{\p  } {\p z} $ defines  the characteristic coorientation of 
the fold $\Sigma$, and  the vector field $\frac {\p  }{\p x}\in TV|_\Sigma$ 
defines the index coorientation.

 \mn  
 The normal form \eqref{eq:cusp} can be extended to a neighborhood of $V$ using 
the parametric Morse lemma. However, we will not need it for our purposes. 

  \bp\label{lm:hol-prep} 
  If $\Phi$ is holonomic    then for appropriate 
auxilliary choices the  virtual vertical bundle ${}^\Phi\Vert$ coincides with 
$\Vert$ and the  isomorphism   
$$\Delta_\Phi: \Vert\to{}^\Phi\Vert=\Vert$$
is the identity.     \ep
 \begin{figure}[hi]
\centerline{\psfig{figure=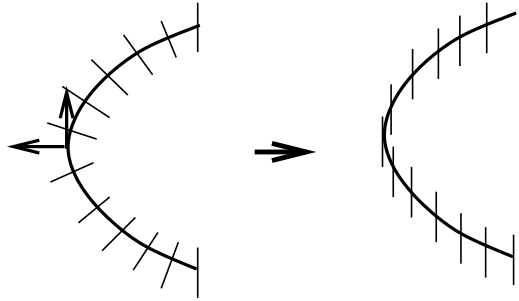,height=40mm}}
\caption {\small Holonomic case: the bundle  ${}^\Phi\Vert$ coincides with 
$\Vert$}
\label{nwi8}
\end{figure}

  {\bf Proof.} Let $\Phi$ be holonomic and $\Phi^0=\phi$.  The bundle $\Vert$ is 
transverse to $V$ over $V\setminus U$, and over $U$  it  splits as 
$\wt\Ver\oplus\wt\lambda$.
We have $\wt\Nor\cap TV=\{0\}$,  the bundle $\wt\lambda$ is tangent to $V$ along 
$\Sigma$ and $\lambda=\wt\lambda|_\Sigma$ is transverse to $\Sigma\,$. 
Let us choose a metric such that   the transversality condition  for the bundles 
$\Vert|_{V\setminus U},\; \wt\Ver|_U,\; \lambda|_U$ are replaced by the 
orthogonality one.
Then  the operator  $\Gamma_\Phi^{-1}$, and hence $\Delta_\Phi$ leaves invariant 
the  bundles $\Vert|_{V\setminus U}$ and $\wt\Ver|_U$. Moreover, 
 on  both these  bundles the operators 
$\Phi^2=d^2\phi$  and  $\Gamma_\Phi $    coincide, and hence 
$\Delta_\Phi=\Id$.

It remains to analyze  $\Delta_\Phi|_{\wt\lambda^+}$.  By definition,    
$$\Delta_\Phi(\wt\lambda^+(v))=\wt 
c(v)\left(\cos\theta(v)\tau^+(v)+\sin\theta(v)\wt n^+(v)\right),\; v\in U,$$ where 
$\wt n^+ =\Gamma_\Phi^{-1}(\wt\lambda^+)$. It is sufficient to ensure
 that the  line  $\Delta_\Phi|_{\wt\lambda}$ coincides with $\wt\lambda^+$ 
because then the  similar equality  for vectors could  be achieved just by 
choosing an appropriate  amplitude function   $\wt c$ in the definition of the 
operator $\Delta_\Phi$.
Note that  we have $\Delta_\Phi(\wt\lambda(v))=\lambda(v)$ for $v\in \p U$ or 
$v\in\Sigma$. To ensure this equality on the rest of $U$ we need to further 
specify   our choices.
As it was explained above in Section \ref{sec:normal-form} we can assume that  
the function $\phi$ in a neighborhood $\Omega\supset U$ in $W$ is given by the 
normal form \eqref{eq:cusp}.  Choosing $\Omega=\{|x|,|z|\leq\eps\}$ we have 
$$U:=V\cap\Omega=\{z=x^2,y_\pm=0, |x|\leq\eps\}\,,$$
and bundles  $\Vert$, $\wt\Ver$ and $\wt\lambda$ are given,
 respectively, by restriction to $V$ of the projections      
$(\sigma,x,z,y_+,y_-)\mapsto(\sigma,z)\,$,\; 
$(\sigma,x,z,y_+,y_-)\mapsto(\sigma,x,z)\,$\;
and $(\sigma,x,z,y_+,y_-)\mapsto(\sigma,z,y_-,y_+)\,.$
 Let us choose the tangent to $V$ vector field
  $  \frac{\p}
 {\p x}+2z\frac{\p}{\p z} $  as $\wt\tau^+$ and recall that we have chosen a 
metric for which the vectors $\tau^+(v)$  and $\wt \lambda^+(v)$ for $v\in\p U$  
are orthogonal.  Let us choose any vector field 
$\wh\nu\in P:=\Span(\frac{\p}{\p x},\frac{\p}{\p z})$ such that 
 \begin{itemize}
\item  $\wh \nu^+|_{\p U_+}=\wt\lambda^+|_{\p U_+}$;
\item $\wh \nu^+|_{\p U_-}=-\wt\lambda^+|_{\p U_-}$;
\item $\wh \nu^+|_{\Sigma}=-\frac{\p}{\p z}$ defines the characteristic 
coorientation;
\item the vector field $\lambda^+|_{\Int U^+}$ belongs to the positive cone 
generated by   $\wt\tau^+$ and $\wh\nu^+$;
\item the vector field $\lambda^+|_{\Int U^-}$ belongs to the positive cone 
generated by   $\wt\tau^+$ and $-\wh\nu^+$.
 \end{itemize}
 Let us pick a metric on $P$ for which  the vector fields  $\wt\tau^+$ and 
$\wh\nu^+$ 
are orthogonal and the vector fields  $\wt\tau^+$ and $\wt\lambda^+$  have length 
1.
 By rescaling, if necessary,  the vector field  $\wh\nu^+$ we can arrange  that 
it has length 1 as well.
  Let us denote by $\theta(v)$ the angle between the vectors $\tau^+$ and 
$\lambda^+$ in this metric. If we construct the virtual vertical bundle 
${}^\Phi\Vert$ 
  with this choice of the metric  and the  angle function $\theta$, then the 
condition 
$\Delta_\Phi(\wt\lambda)=\wt\lambda$ will be satisfied.   \qed
 
\bigskip
In all our results below concerning an extension of  a holonomic FLIF from a 
neighborhood of a closed set $A$ we will always assume that over $\Op A$ all the 
necessary special choices are made to ensure the conclusion of Lemma 
\ref{lm:hol-prep}:  the  virtual vertical bundle ${}^\Phi\Vert$ coincides with 
$\Vert$ and the  isomorphism  $\Delta_\Phi: \Vert\to{}^\Phi\Vert=\Vert$   is   the 
identity,  and hence,  according to Lemma \ref{lm:2-bal}, we have
$\Gamma_\Phi=\Pi_{\wh\Phi}$, where $\wh\Phi=\Phi^\Norm$.
 
  \subsection{Balanced  and well balanced FLIFs}\label{sec:prep} 
   
\mn We call a  FLIF $\Phi$ {\it balanced}   if the compositions  
$$\Norm\mathop{\longrightarrow}\limits^{\Gamma_{\Phi}}\Vert\hookrightarrow TW|_V\;\;\;
\hbox{and}\;\;\;
\Norm\mathop{\longrightarrow}\limits^{\Pi_{\wh\Phi}}\,{}^\Phi\Vert \hookrightarrow TW|_V$$
are homotopic in the space of injective homomorphisms $\Vert\to TW|_V$. Here we 
denote by $\wh\Phi$ the FLIF $\Phi^\Norm$. If $\Phi$ is holonomic over $\Op 
A\subset W$ then we say that $\Phi$ is {\it balanced relative} $A$ if  the 
homotopy can be made fixed over $A$. 
 
\mn  Lemma \ref{lm:2-bal} shows that the balancing condition is equivalent to 
the requirement that the composition 
 $\Vert\mathop{\longrightarrow}\limits^{\Delta_\Phi}\;{}^\Phi\Vert\hookrightarrow TW|_V$ is 
homotopic to the inclusion $\Vert \hookrightarrow TW|_V$  in the space of 
injective homomorphisms $\Vert\to TW|_V$. 
 
\mn  Lemma \ref{lm:hol-bal2} shows that  a holonomic   $\Phi$ is balanced. 
Moreover, it is balanced relative to any closed subset $A\subset W$.

\mn We say that a  FLIF  $\Phi$ is {\it well balanced} if it is prepared and  
the isomorphisms  $\Pi_\Phi,\Gamma_\Phi:\Norm\to\Vert$ are homotopic as   
isomorphisms.
  Similarly we define the notion of a FLIF well balanced relative to a closed 
subset $A$.
  
   It is not immediately clear from the definition that
a well balanced FLIF is balanced. The next lemma shows that this is still the 
case.
\bp
A well balanced FLIF is balanced.
\ep

{\bf Proof.}  We need to check  that over $V\setminus \Sigma$ we have   
$\Pi_{\Phi}=\Phi^2\circ  K^{-1}$ and  $\Pi_{\wh\Phi}=\Phi^2\circ \wh K^{-1}$, 
where $K$ is the projection $\Norm\to\Vert$ and $\wh K$ is the projection 
$\Norm\to{}^\Phi\Vert$. We have $\wh K=T\circ K$, where $T:\Vert\to{}^\Phi\Vert$ is 
the projection along $TV$. Hence, we have  
 $\Pi_{\wh\Phi}=\Pi_{\Phi}\circ T$  which implies, in particular, that the 
projection $T$ is non-degenerate over the whole $V$. Hence, the  composition of 
the projection operator $T$ with the inclusion
 ${}^\Phi\Vert\mathop{\hookrightarrow }\limits^i TW|_V$ is homotopic to the  
inclusion  
  $\Vert\mathop{\hookrightarrow }\limits^j TW|_V$ as injective homomorphisms, 
and so do the compositions $i\circ\Pi_{\wh\Phi}$ and $j\circ\Pi_\Phi$.
  \qed

\mn Note that for the codimension 1 case, i.e. when $n=1$ the well balanced
 condition for a prepared FLIF is very simple:
\bp\label{lm:1-bal}{\bf (Well-balancing criterion in codimension 1)} 
Suppose $\dim\,\zeta=1$. Then any  prepared $\zeta$-FLIF $ \Phi$ 
is well balanced if and only if at 
one point $v\in V\setminus\Sigma$ of every connected component of  $V$   the map 
 $$ (\Pi_\Phi)_v\circ(\Gamma_\Phi)_v^{-1} :\Vert_v \to\Vert_v$$
    is a multiplication by a positive number.
 The same statement holds also  in the relative case.
\ep

\bp\label{lm:forward-balanced}{\bf (Well balanced FLIFs and folded isotopy)}  Let $\Phi$ be a well balanced FLIF. Let 
$h_s:W\to W$ be a diffeotopy,
   $\zeta_s$   a family of $n$-dimensional distributions on $W$, and 
$\,\Theta_s:\zeta_0\to\zeta_s$ a family of bundle isomorphisms covering $h_s$, 
$s\in[0,1]$, such that $h_0=\Id$ and for each $s\in[0,1]$
   \begin{itemize}
    \item submanifold $V_s:=h_s(V) \subset W$ is folded with respect to 
$\zeta_s$ along  $\Sigma_s:=h_s(\Sigma)$;
    \item $dh_s(\zeta_0\cap TV))= dh_s(\zeta_s)\cap TV_s$;
    \item $dh_s|_{\zeta_0\cap TV}=\Theta_s|_{\zeta_0\cap TV}$.
    \end{itemize}
Then the  push-forward  $\zeta_s$-FLIF  $\Phi_s:= (h_s,\Theta_s)_*\Phi$ ,
$s\in[0,1]$, is  
 well balanced.
\ep
{\bf Proof.}
By assumption   $V(\Phi_s)$  is folded with respect to $\zeta_s$.  Next, we observe that all co-orientations  cannot    change in the process of a continuous deformation, and similarly,
  the isomorphisms $\Pi_{\Phi_s} $  and $\Gamma_{\Phi_s}$ vary continuously, 
and hence   remain homotopic as bundle isomorphisms
  $\Norm(\Phi_s)\to \Vert(\Phi_s)$. Thus the well balancing condition is preserved.
  \qed

\mn Note that if $\Phi$ is balanced then the homomorphism 
$\Delta_\Phi:\zeta|_V\to {}^\Phi\Vert $ composed with the inclusion 
${}^\Phi\Vert \hookrightarrow TW$ extends to an injective homomorphism $F:\zeta\to 
TW$. Then
$(\Id,F)_*\Phi$ is a $\nu$-FLIF extending the local $\nu$-FLIF $\wh \Phi$. Here 
we denoted by $\nu:=F(\zeta)$.
 

\bp\label{lm:norm-prep} The $\nu$-FLIF $\,\wh\Phi=\Phi^\Norm$  on $\Op V$ is well 
balanced.  \ep

{\bf Proof.}  We already proved in \ref{lm:tang-to-twist}  that $\wh\Phi$ is 
prepared.  Let us show that $\Pi_{\wh\Phi}=\Gamma_{\wh\Phi}$. According to the 
definition of the push-forward operator we have
$\Gamma_{\wh\Phi}=\Delta_\Phi\circ\Gamma_\Phi$. But according to Lemma 
\ref{lm:2-bal} we have $ \Delta_\Phi\circ\Gamma_\Phi=\Pi_{\wh\Phi}$.  
\qed


\mn Consider a $\zeta$-FLIF $\Phi$. Suppose there exists  a $(k+1)$-dimensional
 submanifold $Y\subset W$,  $Y\supset V$, such that
\begin{itemize}
\item  $Y$ is transverse to $\zeta$;
\item the line field $\mu|_V\subset\Vert$ is an 
eigenspace field for $\Phi^2$, where we denoted $\mu:=\zeta\cap TY$;
\item $\Phi^2|_{N:=\mu^\perp|_V}$ is non-degenerate, where $\mu^\perp$
 is the orthogonal complement to $\mu$ in $\zeta|_Y$.
\end{itemize}
Consider the restriction $\mu$-FLIF $\wt\Phi=\Phi|_Y$ defined as follows:
$\wt\Phi^0=\Phi^0|_Y$, $\wt\Phi^1$ is the projection of $\Phi^1$ along 
 $\mu^\perp$,\; $\wt \Phi^2=\Phi^2|_{\mu}$\,, $\wt\lambda=\lambda$. Note that 
 we have $V(\wt\Phi)=V$ and $\Sigma(\wt\Phi)=\Sigma$.

\mn
We will assume that the bundle $N$ is orthogonal to $TY$. Under this
 assumption we have  $\Gamma_\Phi(N)=N$.
The next criterion for a FLIF to be well-balanced is immediate from the definition.
 
\bp\label{lm:restriction}
If $\wt\Phi$ is prepared then so is $\Phi$. If $\wt\Phi$ is well
 balanced and $\Phi^2|_N=\Gamma_\Phi|_N$ then $\Phi $ is well balanced as well.
\ep 
 

\subsection{Pleating a  FLIF}\label{sec:pleat-FLIF}
We adjust in this section the pleating construction defined in
 Section 
\ref{sec:pleating} for submanifolds  to make it  applicable  for framed well balanced 
FLIFs.
et $ \Phi$ be a well balanced $\zeta$-FLIF.  We will use here the following  
notation from Section \ref{sec:pleating}:
\begin{description}
\item{-\;} $S\subset V_i\subset V\setminus\Sigma, i=0,\dots, n,$ 
is a  closed cooriented codimension $1$ submanifold;
\item{-\;} $U=S\times[-\delta,\delta]\supset S=S\times 0$ is  a tubular 
$\delta$-neighborhood of $S$ in $V_i$;
\item{-\;} $\nu^+\in\zeta$ is  a  unit vector field  defined over a neighborhood  
$\Omega$ of $U$ in $W$; 
\item {-\;} $g:S\times[-\delta,\delta]\times[-\eps,\eps]\to 
\Omega\hookrightarrow W$ is  an embedding such that 
   $\frac{\p g}{\p u}(s,t,u)=
\nu^+(g(s,t,u)),\; (s,t,u)\in S\times[-\delta,\delta]\times [-\eps,\eps]$,
 which maps $S\times 0\times 0$ onto $S$ and  $S\times[-\delta,\delta]\times0$ 
onto $U$;
    \item{-\;}  $\Gamma\subset  P:=[-1,1]\times[-1,1]$ is  an embedded connected 
curve which near 
$\p P$ coincides with  the line $\{u=0\}$;
\item{-\;}  $\wt V\subset W$  is  the result of $\Gamma$-pleating of $V$      
over $S$ in the direction of  the vector field $\nu_+$.
\end{description}
 
 \mn We will make the following additional assumptions:  
 \begin{description}
 \item{$\ast$}
  the  splitting  $\Vert|_S=\Vert_+|_S\oplus \Vert_-|_S$ is extended  to  a 
splitting $\zeta=\zeta_+\oplus\zeta_-$ over the neighborhood
 $\Omega\subset W$;
\item{$\ast$}  the vector field  $\nu^+$ is   a section of  either $\zeta_-
|_\Omega$ or $\zeta_+|_\Omega$;
 \item{$\ast$}  the vector field $\nu_+|_U$ is an eigenvector field for  
$\Phi^2$;
 \item{$\ast$} $\Norm(\Phi)|_U=\Vert(\Phi)|_U$ and  
$\Delta_\Phi|_{\Vert|_U}=\Id$.
  \end{description}
  
  There exists a diffeotopy $h_s:W\to W$   supported in   $\Omega$ connecting 
$\Id$ with a diffeomorphism $h$    such that $h(V)=\wt V$. We denote $\wt U=\wt U_\Gamma:=h_1(\Gamma)$. Let 
  $\Psi_s:\zeta\to\zeta$, $s\in[0,1]$, be a family of isomorphisms covering
$h_s$ which preserve  $\Vert_\pm$ and $\nu^+$. 

\mn The manifold $\wt U$ is folded with  respect to $\zeta$ with the 
fold   $\wt S=\bigcup\limits_1^{2N}\wt S_j$ where $\wt S_{j}= h_1(S_j)$, where 
$S_j=S\times t_j$,  $-\delta<t_1<\dots t_{2N}<\delta$.  Over $\wt S$ we have 
$\wt\tau=\nu= T\wh V\cap \zeta$. 
  
 \mn Consider the push-forward FLIF  $\overline\Phi:=(h_1,\Psi_1)_*\Phi$.
Though the manifold $V(\overline\Phi)=\wt V$ is folded with respect to $\zeta$,
it is not prepared. We will  modify 
$\overline\Phi$ to a prepared FLIF $\wt\Phi=\Pleat_{S,\nu^+,\Gamma}(\Phi)$ as follows.

 \mn Let $\wt c:\wt U \to \bbR$ be  a function which on $\p \wt U=\p U$
 coincides with the eigenvalue function   of the operator $\Phi^2$ for  
 the eigenvector field $\nu^+$, and have the fold   
 $\wt S:=\bigcup\limits_1^{2N}\wt S_j$ 
  as its regular $0$-level. We call component of $\wt U\setminus \wt S$
  {\it  positive} or {\it negative} depending on the sign of the function $\wt c$ on
    this component.  We then define
 \begin{figure}[hi]
\centerline{\psfig{figure=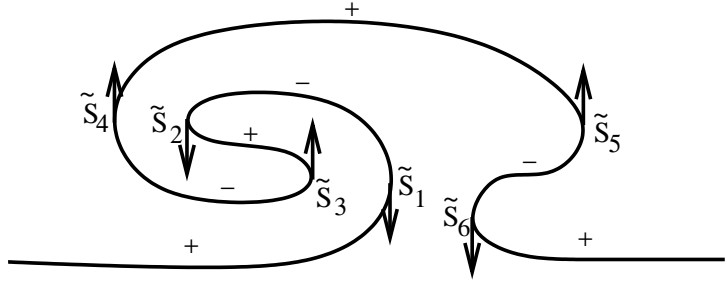,height=40mm}}
\caption {\small $\Gamma$-pleating of a well balanced FLIF }
\label{nwi10}
\end{figure}  
\begin{itemize}
\item $ \wt\Phi^1=\overline\Phi^1$;
\item $\wt\Phi^2|_{\nu^\perp}=\overline\Phi^2|_{\nu^\perp}$;
\item $\wt\Phi^2(\nu^+)=\wt c\,\nu^+$;
\item $\lambda^+(\wt\Phi^2)=\pm\nu^+$,
where the sign is chosen in such way that the vector field $\lambda^+(\wt\Phi^2)$
define an inward coorientation of positive components of $\wt U\setminus \wt S$, 
see Fig.  \ref{nwi10}.
 \end{itemize}  

We say that $\wt\Phi=\Pleat_{S,\nu^+,\Gamma}(\Phi)$ is obtained from $\Phi$ by {\it  $\Gamma$-pleating over
 $S$ in the direction of the vector 
field $\nu^+$}
  see Fig. \ref{nwi10}.

\bp\label{lm:pleat-bal}   
 The FLIF  $\wt\Phi$ is well balanced. 
\ep 
{\bf Proof.}    
 Consider the $(k+1)$-dimensional manifold 
 $$Y:=g(S\times[-\delta,\delta]\times[-\eps,\eps])\subset W\,.$$ 
 Then $Y$ is transverse $\zeta$ and $\zeta\cap TY=\nu$.  We also note that 
 the orthogonal complement $\nu^\perp$ of $\nu\in\zeta$ is orthogonal to $TY$,
  $\wt \Phi^2|_{\nu^\perp}=\Pi_{\wt\Phi}|_{\nu^\perp}=\Gamma_{\wt\Phi}|_{\nu^\perp}$.
   According  to \ref{lm:restriction} it is sufficient to check that the restriction
   $\wh \Phi:=\wt\Phi|_Y$ is well balanced, rel. $\p Y$.
  First, we need to check that this restriction is prepared. By construction,
   $\wt V=V(\wh\Phi)$ is folded with respect to $\nu$ and the vector field 
   $\lambda^+(\wh\Phi)=\lambda^+(\wt\Phi)$ defines the index coorientation 
   of $\wt S$ in $\wt V$.
   Next, we need to check that the vector field
    $n^+(\wh\Phi)=n^+(\wt\Phi)=\Gamma_{\wt\Phi}^{-1}(\lambda^+(\wt \Phi)$
     defines the characteristic coorientation of the fold.  It is sufficient
      to consider the case when  $S$ is the point, and hence $\dim \,Y=2$. 
       The general picture is then obtained by taking a direct product with $S$.
    \begin{figure}[hi]
\centerline{\psfig{figure=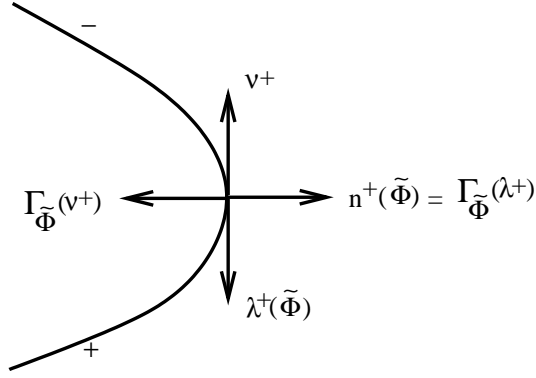,height=50mm}}
\caption {\small Vector $n^+(\wt\Phi)$ determines the characteristic
 coorientation of the fold}
\label{nwi12}
\end{figure} 
Note that the characteristic co-orientation of the fold $\wt S_j$ is given
 by the vector field $\frac{\p}{\p t}$  if $j$ is odd,
 and by  $-\frac{\p}{\p t}$ if $j$ is even.  Consider first the case when $j$ is odd,
 see Fig. \ref{nwi12}.
Then  if the lower branch of the parabola is positive then the vector field  
 $\Gamma_{\wt\Phi}(\nu^+)$ defines the same coorientation as the vector field 
 $-\frac{\p}{\p t}$. But in this case $\lambda^+=-\nu^+$, and hence
 $\Gamma_{\wt\Phi}(\lambda^+)$ defines the characteristic coorientation of the fold.
  The other cases can be considered in a similar way.  
Finally, we use Lemma   
     \ref{lm:1-bal} to  conclude that
    $\wh\Phi$ is well balanced relative the boundary $\p Y$.
  \qed
  
\begin{figure}[hi]
\centerline{\psfig{figure=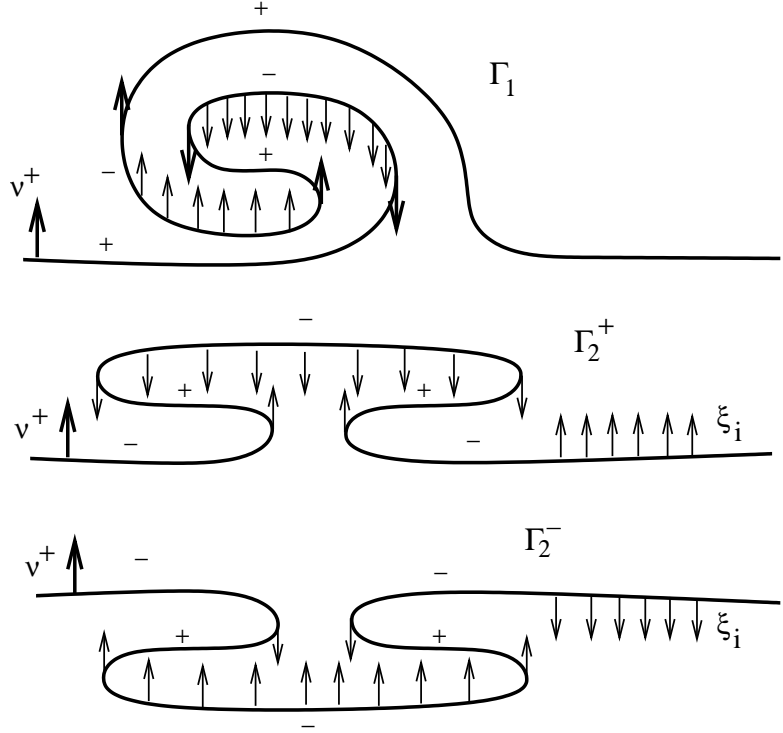,height=100mm}}
\caption {\small  Framing of a $\Gamma$-pleated FLIF }
\label{nwi11}
\end{figure}

\mn In order to extend the $\Gamma$-pleating operation to framed well-balanced 
FLIFs we need to impose additional constraints on the choice of the vector field 
$\nu^+$ and the curve $\Gamma$, see Fig. \ref{nwi11}.  For each $j=1,\dots 2N$
 denote by $\sigma_j$ the proportionality coefficient in
 $\lambda^+|_{\wt S_j}=\sigma_j\nu^+|_{\wt S_j}$, $\sigma_j=\pm 1$.
  Then we require that
 \begin{description}
 \item{($\alpha$)}  if $\wt S_j$ and $\wt S_{j+1}$, $j=1,\dots, 2N-1$ bound
  a negative component of  $\wt U\setminus\wt S$  then $\sigma_j=\sigma_{j+1}$;
 \item{($\beta$)} if the component bounded by $\wt S_1$ and $\wt S_{2}$ is positive
 then $\sigma_1=\sigma_{2N}=\pm1$ for $\nu^+=\pm\xi^i$.
 \end{description}

\bp\label{lm:framed-pleat} {\bf (Pleating a framed FLIF)} If  $\,\nu^+$ and 
$\Gamma$  satisfy the above conditions, then given a framed 
well balanced FLIF $(\Phi,\xi)$ the FLIF $\wt\Phi=\Pleat_{S,\nu^+,\Gamma}$ 
admits a framing $\,\wt\xi$, where the framing $\wt\xi$ coincides with $\xi$ outside $\wt U$.
\ep
 
{\bf Proof.} The proof is illustrated on Fig.~\ref{nwi11}.
 In the case  $\nu^+\in\Vert_+$ the pleating construction adds a $1$-dimensional 
negative eigenspace to $\Vert_-$ restricted to negative components
  $\wt U\setminus\wt S$.   Condition $(\alpha)$ then 
allows us to frame this 1-dimensional space either with $\xi^{i+1}:=\sigma_j\nu^+$.
Similarly, if $\nu^+\in\Vert_-$ (or, equivalently when the component bounded
 by $\wt S_1$ and $\wt S_2$ is positive)
then the  pleating construction removes  the 
negative eigenspace generated by $\nu^+$ on positive components.
The  remaining negative components bounded $\wt S_{2j}$ and $\wt S_{2j+1}$, $j=1,\dots, N-1$, can be framed with  $\wt\xi:=(\xi_1,\dots,\sigma_{2j}\xi_i)$.   
Condition {($\beta$)} ensures that the existing framing in the complement of $U$ 
satisfies the necessary boundary conditions on $\wt S_1$ and $\wt S_{2N}$.
\qed
 
\mn \bp\label{lm:standard-curves} Given any framed well balanced FLIF 
$(\Phi,\xi)$, one of  the curves $\Gamma_1,\Gamma_2^\pm$ shown on Fig. 
\ref{nwi5} 
can always be used as the curve $\Gamma$ to produce a framed well balanced 
FLIF $(\wt\Phi,\wt\xi)$ by a $\Gamma$-pleating.  
\ep
{\bf Proof.} Indeed, as it follows from Criterion \ref{lm:framed-pleat},
the curve $\Gamma_1$ can   always be used if $\nu^+|_S\in\Vert_+$, while if 
$\nu^+|_S\in\Vert_-$ then the curve  $\Gamma_2^\pm$  can be  used 
in the case $\nu^+=\pm\xi^i$, see 
Fig.~\ref{nwi11}.
\qed

\mn
The next proposition is a corollary of Theorem \ref{thm:folds} and the results 
discussed  in the current section.
 \bp\label{thm:folds2}{\bf (Pleated isotopy of framed well balanced FLIFs)}
 Let $\zeta_s$, $s\in[0,1]$, be a family of $n$-dimensional distributions on 
$W$, and 
   $(\Phi,\xi)$   a  framed well-balanced $\zeta_0$-FLIF with $V(\Phi)=V\subset 
W$. Then there exist
\begin{itemize}
\item   a framed well balanced $\zeta_0$-$FLIF$ $\wt\Phi$ obtained from $\Phi$  
by a sequence of pleatings,  and
\item  a $C^0$-small  isotopy $h_s:V\to W$, $s\in[0,1]$ such that $h_0$  is the 
inclusion $V\hookrightarrow W$ and $\wt V_s:=h_s(V(\wt\Phi))$ is folded with 
respect to $\zeta_s$ along $\wt\Sigma_s:=h_s(\Sigma(\wt\Phi))$.
\end{itemize}
If $\Phi$ is holonomic over $\Op A$ then one can arrange  that $\wt\Phi=\Phi$ on 
$\Op A$ and that the homotopy $h_s$ is fixed over $\Op A$.
\ep
{\bf Proof.}
According to Theorem \ref{thm:folds} there exists a  manifold  $\wt V$ for which 
the isotopy with the required properties does exist. This manifold  can  be 
constructed beginning  from $V$ by a sequence of  $\Gamma^+_0$-pleatings 
along the boundaries of balls embedded into $V\setminus\Sigma$, in the direction 
of vector fields which extend to these balls. The latter property allows us to 
deform these vector fields into vector fields contained  in $\Vert_+$ or 
$\Vert_-$ (we need  to use $\Vert_-$   only if   $\dim\,\Vert_+=0$). Moreover, when using 
$\nu^+\in\Vert_-|_{V_i}$  and when $i=\dim\,\Vert_-|_{V_i}>1$ we can deform it 
further into the last  vector $\xi^i$ of the framing.  In the case  $i=1$ we can 
deform
$\nu^+$ into  $\pm \xi^i$, but we cannot, in general, control the sign.  Note that we need to use this case only if $n=1$.  
As it  was explained in Remark \ref{rem:special-Gamma},
 we can replace  at our choice 
 each 
$\Gamma^+_0$-pleating  in the 
statement of Theorem \ref{thm:folds}
 by any of  the $\Gamma$-pleatings with 
$\Gamma=\Gamma_1,\Gamma_2^\pm$. 
But according to Lemma
 \ref{lm:standard-curves}  one can always use  one of these curves to pleat in the class of   framed   well balanced   FLIFs.   It remains to observe that  if $\Phi$ is holonomic 
over $\Op A$ then all the constructions which we used  in the proof can be made relative  to $\Op A$. 
\qed

\subsection{Stabilization}\label{sec:stab}

 \mn   Let $ \Phi $ be a  $\zeta$-FLIF.
  Suppose that we are given a connected 
domain $U\subset V\setminus\Sigma$   with smooth boundary  such that the bundles 
$\Vert_\pm|_U$ are trivial. Let    $C$ be an exterior collar  of $\p U\subset 
V\setminus\Sigma$. We set $U':=U\cup C$.     
 
Let us assume that $U$ is contained in $V^i$. If  $i<n$ we choose a section 
$\theta^+$ of the bundle $\Vert_+$ over $U'$ and   we define  a {\it negative   
stabilization}  of $\Phi$ over $U$ as a FLIF $ \wt\Phi=\St^-_{U,\theta^+} 
(\Phi)$ such that 
 \begin{itemize}
 \item $\wt\Phi^1=\Phi^1$;
 \item $\wt\Phi^2=\Phi^2$ over $V\setminus U'$;
 \item $\Sigma(\wt\Phi)=\Sigma(\Phi)\cup\p U$; $\Int U\subset V^{i+1}(\wt\Psi)$;
  \item   $\Vert_-(\wt\Phi)|_{\Int U}=\Span(\Vert_-(\Phi)|_{\Int U},\theta^+)$;
  \item $\lambda^+(\wt\Phi)|_{\p U}=\theta^+$.
   \end{itemize}
   We will omit a reference to $\theta$ in the notation and write simply $\St^-
_{U} (\Phi)$ when this choice will be irrelevant.
   
  \mn Note that in order to construct $\wt\Phi^2$  on $U'$ which ensures these 
property we need  to adjust the background metric on $\zeta$ to make $\theta^+$ 
an eigenvector field for $\Phi^2$ corresponding the eigenvalue $+1$.
    The vector field $\theta^+$ remains the eigenvector field for $\wt\Phi^2$ 
but the eigenvalue function is changed to $c:U'\to[-1,1]$, where $c$ is negative 
on $U$, equal to $1$ near $\p U'$ and has $\p U$ as its regular $0$-level.  
   
 \mn  If    the FLIF $\Phi$  is framed by   $\xi=(\xi^1,\dots,\xi^i)$ then 
$\wt\Phi$ 
can be canonically framed by $\wt\xi$ such that $\wt\xi=\xi$ over $V\setminus U$ 
and $\Vert_-(\wt\Phi)|_{\Int U}$ is framed by $\wt\xi:=(\xi^1.\dots,\xi^i,\theta^+)$ and  we define 
$$\St^-_U(\Phi,\xi):=(\St^-_{U,\,\xi_i}(\Phi),\wt\xi),$$ 
   
   \mn  In the case when $U\subset V_i$ and $i>0$ we can  similarly define a 
{\it positive 
     stabilization}  of $\Phi$ over $U$ as a FLIF $ \wt\Phi=\St^+_{U,\,\theta} 
(\Phi)$, where $\theta$ is a section of $\Vert_-$ over $U'$ such that 
 \begin{itemize}
 \item $\wt\Phi^1=\Phi^1$;
 \item $\wt\Phi^2=\Phi^2$ over $V\setminus U'$;
 \item $\Sigma(\wt\Phi)=\Sigma(\Phi)\cup\p U$; $\Int U\subset V^{i-1}(\wt\Psi)$;
   \item    $\Vert_+(\wt\Phi)|_{\Int U}=\Span(\Vert_+(\Phi)|_\Int U,\theta^+)$;
     \item $\lambda^+(\wt\Phi)|_{\p U}=\theta^+$.
   \end{itemize}
   If $\Phi$ is framed by a framing $\xi=(\xi^1,\dots,\xi^i)$ then we will 
always choose $\theta^+=\xi^i|_U$ 
   and define a  positive stabilization by the formula
     $$\St^+_U(\Phi,\xi):=(\St^+_{U,\,\xi_i}(\Phi),\wt\xi),$$ 
     where $\wt\xi|_{\Int U}=(\xi^1,\dots,\xi^{i-1})$.

  \bp\label{lm:stab-to-balanced} {\bf (Balancing via stabilization)}    Any  FLIF can be stabilized to a  balanced 
one. If $ \Phi$ is balanced and   $\chi(U)=0$ then $\St^\pm_U(\Phi)$ is balanced 
as well. The statement holds also in the relative form. \ep
  {\bf Proof.}
The obstruction  for  existence of a fixed over $A\subset V$ homotopy between  
two monomorphisms $\Psi_1,\Psi_2:\Norm\to TW|_V$ is an $n$-dimensional 
cohomology 
class   $  \delta(\Psi_1,\Psi_2;V,A)\in H^k(V,A;\pi_k(V_n(\bbR^{n+k})))$, or 
more precisely a cohomology class with coefficients in the local system 
$\pi_k(V_n(T_vW)), v\in V$. Note that $\pi_k(V_n(\bbR^{n+k}))=\bbZ$ if $k$ is 
even or $n=1$ and $\bbZ/2$ otherwise.  It is straightforward to see that 
$$\delta(\Delta(\Phi),
\Delta(\St^\pm_U(\Phi);U,\p U)=\begin{cases}\chi(U)\Theta,& k\;\hbox{is 
even};\cr
\pm\chi(U)\Theta,& k\;\hbox{is odd},
\end{cases}
$$ for an appropriate choice of  a generator $\Theta$   of  $H^k(U,\p 
U;\pi_k(V_n(\bbR^{n+k})))$.
  Hence, stabilization over a domain with vanishing Euler characteristic does 
not change the obstruction class $\delta(\Gamma(\Phi),\Delta(\Phi))$ and  with 
the exception of the case $k= n=1 $   this obstruction class  can be changed in 
an arbitrary way by an appropriate choice of $U$.  Indeed, if $k>1$ then one can 
take as $U$ either the union of $l$ copies of $n$-balls or a  regular 
neighborhood of an embedded bouquet  of $l$ circles (comp. a similar argument in 
\cite{[EGM11]}). If $k=1$ and $n>1$ then the sign issue is irrelevant because 
the 
obstruction is 
  $\bbZ/2$-valued.
  If $k=n=1$ then one may need two successive stabilizations in order to balance 
a FLIF.
 Indeed, the domain $U$ in this case is a union of some number $l$ of intervals, 
and hence $\chi(U)=l$. Thus the positive stabilization increases the obstruction 
class by $l$, while the negative  one decreases it by $l$. Suppose, for 
determinacy, we want to stabilize over a domain in $V_0$.  If we need to  change  
the obstruction class by $-l$ then we just negatively stabilize over the union 
of $l$ intervals. If we need to change it  by  $+l$ we first negatively  
stabilize 
over one interval $I$ and then positively stabilize over the union of $l+1$ 
disjoint intervals in $I$. \qed
 
 \subsection{From balanced to  well balanced FLIFs}\label{sec:bal-well}
 \bp\label{prop:main}{\bf (From balanced to well balanced)}
   Let $(\Phi,\xi)$ be a balanced framed $\zeta$-FLIF which 
is holonomic  over a neighborhood of a closed  subset $A\subset W$. Then  there 
exists a framed  well-balanced FLIF $(\Phi',\xi')$ which coincides with $\Phi$ 
over $\Op A$. In addition, $V(\Phi')$ is obtained from $V(\Phi)$ via a $C^0$-
small, fixed on $\Op A$ isotopy. \ep 
 
  {\bf Proof.}  
    There exists  a family    of monomorphisms $\Psi_s:\Vert\to TW$, 
$s\in[0,1]$, connecting 
  $\Vert\mathop{\longrightarrow}\limits^{\Delta_\Phi}
   {}^\Phi\Vert\hookrightarrow\tau$
and the inclusion $j:\Vert\hookrightarrow\tau$.  
The homotopy can be chosen fixed over $\Op A$.
The family $\Psi_s$ can be extended to a family of monomorphisms $\zeta\to TW$. 
We will keep the notation $\Psi_s$ for this extension. Denote 
$\zeta_s:=\Psi_s(\zeta)$, $s\in[0,1]$. Thus $\zeta_1=\zeta$ and $\zeta_0$ is an 
extension to $W$ of the bundle $\Norm^\Phi$.
Lemma \ref{lm:norm-prep} then guarantees that the push-forward $\zeta_0$-FLIF 
$(\Id,\Psi_0)_*(\Phi,\xi)$ is well balanced.
    According to Theorem \ref{thm:folds2} there exists  a  well balanced 
 framed   $\zeta_0$-FLIF $(\wh\Phi,\wh \xi)$ where $\wh V=V(\wh\Phi)$ is 
obtained from $V$ by a $C^0$-small isotopy  which i8s fixed 
outside a neighborhood of $V$ and over a neighborhood of $A$,  and      a $C^0$-
small supported in $(\Op\wh V)\setminus A$ isotopy $g_s $ starting with  
$g_0=\Id$    such that for each $s\in[0,1]$ the manifold  $\wh V_s:=g_s(\wh V)$ 
is folded with respect to   $\zeta_s$ along  $\wh\Sigma_s=g_s(\wh\Sigma)$. There 
exists a family of bundle isomorphisms $\Theta_s:\zeta_0\to\zeta_s$ covering the 
diffeotopy $h_s$ and such that $\Theta_0=\Id$ and $\Theta_s=dg_s$ over  the line 
bundle $TV|_{\wh\Sigma}\cap\zeta_0$.  The homotopy $\Theta_s$ can be chosen 
fixed over $\Op A$.
 Then, according to Lemma  \ref {lm:forward-balanced},  the push-forward  
$\zeta$-FLIF
 $(g_1,\Theta_1)_*(\wh\Phi,\wh\xi)$  is well balanced relative $A$.       
 
   \qed

 \subsection{Formal extension}\label{sec:formal}
 
\bp
\label{thm:formalext2}{\bf (Formal extension theorem)}
Any framed $\zeta$-FLIF  $(\Phi, \xi )$ on $\Op A\subset W$ extends to a  framed $\zeta$-FLIF 
$\,(\wt\Phi,\wt\xi)$ on the whole manifold  $W$.  
\ep 
The proof is essentially   Igusa's argument in  \cite{[Ig87]} (see pp.438-442).

\begin{figure}[hi]
\centerline{\psfig{figure=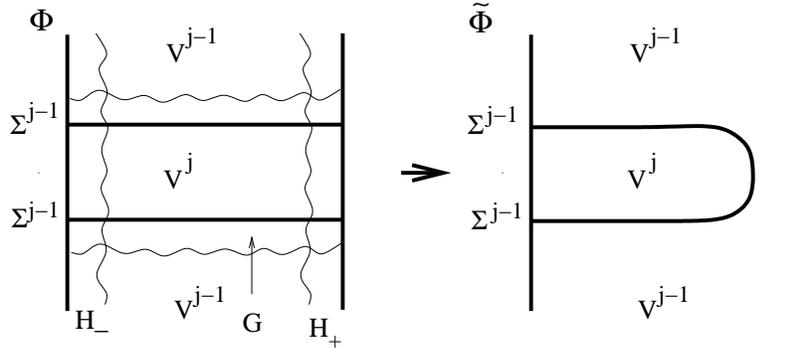,height=50mm}}
\caption {\small  Decreasing the negative index}
\label{nwi13}
\end{figure}

\mn We begin with the following lemma which will be used as an induction step in 
the proof.
\bp\label{lm:formalext} {\bf (Decreasing the negative index)} Let $j=1,\dots, n$. Suppose $W $  is a cobordism 
between $\p_-W$ and $\p_+W$, and for a framed FLIF $(\Phi,\xi)$   on $W$ one has 
$V^i=\varnothing$  for $i> j$. Then there exists a framed FLIF  $(\wt 
\Phi,\wt\xi)$  such  that 
\begin{itemize}
\item $\Phi=\wt\Phi$ on $\Op (\p_-W)$;
\item   $V^i(\wt\Phi)\cap \p_+W=\varnothing $ for $i\geq j$.
\end{itemize}
\ep
{\bf Proof of \ref{lm:formalext}.}    
 To prove the claim we recall that the  $j$-dimensional bundle  $\Vert_-$  over 
$V^j$ is 
 trivialized by the framing   $\xi=(\xi^1,\dots,\xi^j)$, 
 and $\xi^j|_{\Sigma^{j-1}}=\lambda^+$.
  We can extend  the vector field  $\xi^j$ to a 
  neighborhood $G$ of $V^j$ in $V^{j-1}\cup V^j\cup\Sigma^{j-1}$  as a unit 
vector field in $V^{j-1}_+$.
Let $X^j$ be  the self-adjoint linear operator $\Vert_+\to\Vert_+$ defined on 
the neighborhood $G$   which orthogonally  projects $\Vert$ to the line bundle 
spanned by $\xi^j$. Choose neighborhoods $H_-\supset\p_-W$ and $H_+\supset\p_+W$ 
in $W$ with disjoint closures and consider  a cut-off function 
$\theta:V\to\bbR_+$ which is equal to  $0$ on $(V\cap H_-)\cup(V\setminus G)$ 
and  equal to $1$ on $V^j\cap H_+$.
Set  $\wt\Phi^2:=\Phi^2+C\theta X^j$. Then for a sufficiently large $C>0$ the self-
adjoint operator $\wt\Phi^2$  coincides with $\Phi^2$ on $V\cap H_-$, has 
negative index $\leq j$ everywhere, and $<j$ on $V^j\cap H_+$,\, see Fig.\ref{nwi13}.

The kernel of 
$\wt\Phi^2$ on $\Sigma^{j-1}(\wt\Phi)$ is generated by $\xi^j$, and hence there 
is a canonical way to define the vector field $\lambda^+(\wt\Phi)|_{\Sigma^{j-1}}$.
 Note that   $\wt\Phi 
=\Phi$ in the complement $V\setminus V^{j-1}(\Phi)\cup V^j(\Phi)\cup\Sigma^{j-
1}(\Phi)$ and at each point $v\in V^{j-1}(\Phi)\cup V^j(\Phi)\cup\Sigma^{j-
1}(\Phi)$  the negative eigenspace $ \Vert_-(\wt\Phi)$  coincides  either with
$\Vert_-(\Phi)$, or with the span of the vectors $\xi^1,\dots,\xi^{j-1}$. Hence, 
the framing $\xi$ of $\Phi$ determines a framing $\wt\xi$ of $\wt\Phi$.
\qed

 \mn{\bf Proof  of  \ref{thm:formalext2}.}
Let $\Phi =(\Phi^0,\Phi^1,\Phi^2,\lambda^+)$.  
Without loss of generality we can assume that $(\Phi ,\xi )$ is defined on an
 $(n+k)$-dimensional domain $C\subset W, \Int C\supset A$, with smooth boundary.  
Note that if $\Phi^2|_{V\cap\p C}$ is positive definite, then the extension 
obviously exist. Indeed, we can extend $\Phi^1$ in any generic way to $W$, and 
then extend $\Phi^2$ as a positive definite operator on $\Vert$. We will 
inductively reduce the situation to this case.
 Let $C'\subset\Int C$ be a smaller domain such that $A\subset\Int C'$. 
   Let us apply \ref{lm:formalext} to the cobordism $W_{(0)}=C\setminus \Int C'$
    between $\p_-W_{(0)}=\p C'$ and $ \p_+W_{(0)}=\p C$ and to the restriction   
     $ (\Phi,\xi)|_{W_{(0)}}$ in order to  modify   $ (\Phi,\xi)|_{W_{(0)}}$ 
into a framed FLIF $(\Phi_{(0)},\xi_{(0)})$ which 
     coincides with $(\Phi,\xi)  $ near $\p_- W_{(0)}$ and such that 
     $V^n( \Phi_{(0)})\cap (\p_+W_{(0)})=\varnothing$. Then for  a sufficiently 
small
tubular neighborhood $W_{(1)}$ of $\p_+W_{(0)}$ in $W_{(0)}$ we
 have $W_{(1)}\cap A=\varnothing$ and $W_{(1)}\cap V^n(\Phi_{(0)})=\varnothing$. 
 We view $W_{(1)}$ as a cobordism between $\p_-W_{(1)}=\p 
W_{(1)}\setminus\p_+W_{(0)}$ and $\p_+W_{(1)}  =\p_+W_{(0)}$. Now we  again 
apply \ref{lm:formalext}
  to the cobordism $W_{(1)}$ and  $ \Phi_{(0)}|_{W_{(1)}}$ and construct a 
framed FLIF $(\Phi_{(1)},\xi_{(1)})$ on $W_{(1)}$ which coincides  with 
$(\Phi_{(0)},\xi_{(0)})$ near $\p_- W_{(1)}$
  and such that $V^i(\Phi_{(1)})\cap \p_+W_{(1)}=\varnothing)$ for $i\geq n-1$. 
Continuing this process we construct a sequence of nested cobordisms $C\supset 
W_{(0)}\supset W_{(1)}\supset\dots\supset W_{(n-1)}$ and a sequence of framed 
FLIFs $(\Phi_{(j)}, \xi_{(j)})$ on $W_{(j)}$, $j=0,\dots,n-1$,  such that  for 
all $j=0,\dots,n-1 $
  \begin{itemize}
  \item $\p_+W_{(j)}=\p C$;
  \item $(\Phi_{(j+1)},\xi_{(j+1)})$ coincides with $(\Phi_{(j)},\xi_{(j)})$ on 
$\Op(\p_-W_{(j+1)})$;
  \item $V^i(\Phi_{(j)})\cap \p_+W_{(j)}=\varnothing$ for $i\geq n-j$.
  \end{itemize}
  Let us also set $W_{(n)}=\varnothing$.
Hence we can define a   framed formal Igusa function  $(\wt\Phi,\wt\xi)$ over 
$C$  by setting
 $(\wt\Phi,\wt\xi)=(\Phi,\xi)$ on  $C'$ and 
$(\wt\Phi,\wt\xi)=(\Phi_{(j)},\xi_{(j)})$ on $W_{(j)}\setminus W_{(j+1)}$ for 
$j=0,\dots, n-1$. 
 Note that the quadratic part $\wt\Phi^2$ of $\wt\Phi $ is positive definite on 
$\p C$, and hence the framed formal Igusa function 
 $\wt\Phi $ can be extended to the whole $W$.
  \qed
 
 \subsection {Integration near $V$}
 \bp\label{lm:loc-int} {\bf (Local integration of a well balanced FLIF)}
 Any well balanced   $\F$-FLIF $\Phi $ can be made 
holonomic  near $V$ after a small perturbation near $V$. Namely,  there exists 
a  homotopy   of well balanced FLIFs $\Phi_s, s\in[0,1]$, $s\in[0,1]$, beginning 
with $\Phi_0=\Phi$ with the following properties:
 \begin{description}
\item{-} $V( \Phi_s)=V(\Phi), \Sigma(\Phi_s)=\Sigma(\Phi)$ for all $s\in[0,1]$;
\
\item{-} $( \Phi_s^2,\lambda^+_s)$ is $C^0$-close  to $(\Phi^2, \lambda^+)$ for 
all $s\in[0,1]$;
\item{-} $\Phi_1$ is holonomic on $\Op V$.
 \end{description}
 If for a closed subset $A\subset W$ the FLIF $\Phi$ is already holonomic over
  $\Op A\subset W$ then the homotopy can be chosen fixed over $\Op A$.
 \ep
 
 {\bf Proof.} According to Lemma \ref{lm:fold-normal} there exist local 
coordinates
 $(\sigma,t,z,y)$ in a neighborhood of $\Sigma$ in $W$, where
  $\sigma\in\Sigma,\; x,z\in\bbR$ and $y\in\Ver|_\Sigma$ such that  the manifold 
$V$ is given by the equations
$z=x^2,y=0$ and the foliation $\F$  is given by the fibers of the projection    
$(\sigma,x,z,y)\to(\sigma,z)$. 
The vector field $\frac{\p}{\p x}$ generates the line bundle 
$\lambda=TV|_\Sigma\cap\Vert$ and we can additionally arrange that  
$\frac{\p}{\p x}|_\Sigma=\lambda^+$.
By a small $C^0$-small perturbation of the operator $\Phi^2$  (without changing 
it along $\Sigma$) we can arrange the the vector field $\frac{\p}{\p x}$ serves 
an an eigenvector field for $\Phi^2$  in a neighborhood  of $\Sigma$.  We will keep the notation $\lambda$  for the extended line field   $\frac{\p}{\p x}$. Then the operator 
$\Phi^2:\Vert=\Ver \oplus\lambda\to\Ver\oplus\lambda$ can be written as $A\oplus c$, 
where $A$ is a non-degenerate  self-adjoint operator and $c$ is an operator  
acting on the line bundle $\lambda$ by   multiplication  by a 
function $c=c(\sigma,x)$ on $\Op\Sigma\subset V$ such that  for all 
$\sigma\in\Sigma$ we have $c(\sigma,0)=0$, $d(\sigma):=\frac{\p c}{\p 
x}(\sigma,0)>0$.

Define  a function $\phi$ on $\Op\Sigma\supset W$ given by the formula
\begin{equation}
\label{eq:integration}\phi(\sigma,x,z,y)=\frac{d(\sigma)}6(x^3-
3zx)+\frac12\langle Ay,y\rangle.
\end{equation}
Then $V(\phi)=V\cap\Op\Sigma$ and
 the operator $d^2_\F\phi:\Ver \oplus\lambda\to\Ver\oplus\lambda$  is equal to 
$A\oplus\wh c$, where the operator  $\wh c$ acts on $\lambda$ by  multiplication 
by the function
$d(\sigma)x$. Hence the operator functions $d^2\phi$ and $\Phi^2$ coincides with 
the first jet along $\Sigma$, and therefore, one can adjust $\Phi^2$ by a $C^0$-
small homotopy to make $\Phi^2$ equal to $d^2\phi$ over $\Op\Sigma\subset W$.
To extend $\phi$ to a neighborhood $\Op V\subset W$ we observe that the 
neighborhood   of $V$ in $W$ is diffeomorphic to the neighborhood of the zero 
section in the total space of the bundle $\Vert|_{V\setminus U}$.  In   
  the corresponding  coordinates  we define   
  $\phi(v,y):=\frac12\langle \Phi^2(v)y,y\rangle,$ $v\in V,y\in \Vert_v$.
On the boundary of the neighborhood of $\Sigma$ where we already constructed  
another function,  the two functions differ  in terms of order $o(||y||^2)$. 
Hence they can be glued together   without affecting  $d^2_\F\phi$, and thus we get  a leafwise Igusa function $\phi$ with $d^2_\F\phi=\Phi^2$.
It remains to extend $\nabla_\F\phi$ as a non-zero section of the bundle $T\F$ 
to the whole $W$. According to Lemmas \ref{lm:hol-bal2} and  \ref{lm:hol-prep}  we have 
$\Gamma_\phi=\Pi_\phi=d^2_\phi=\Phi^2$.  Then the well balancing condition for $\Phi$ 
implies that $\Gamma_\phi$ is homotopic (rel. $\Op A$) to $\Gamma_\Phi$ as 
isomorphisms $\Norm\to\Vert$. But this implies that there is a homotopy  (rel. 
$\Op A$) of sections $\Phi^1_s:W\to\Vert, s\in[0,1]$, connecting 
$\Phi^1_0=\Phi^1$ and $\Phi^1_1=\nabla_\F\phi$ and such that the zero set 
remains regular and unchanged. \qed

\section{Proof of    Extension Theorem \ref{thm:main} }\label{sec:proof}

\mn{\sc Step 1. Formal extension.} We begin with a leafwise framed Igusa 
function $(\phi_A,\xi_A)$. Using \ref{thm:formalext2} we extend it to a FLIF 
$(\Phi, \xi)$ on $W$. 

All consequent steps are done without changing anything on $\Op A$.

\mn{\sc Step 2. Stabilization.} Using \ref{lm:stab-to-balanced} 
we make $(\Phi,\xi)$ balanced.

\mn{\sc Step 3. From balanced to well balanced.} Using \ref{prop:main} we 
further improve $(\Phi,\xi)$ making it  well balanced.

\mn{\sc Step 4. Local integration near $V$}. Using \ref{lm:loc-int} we deform 
$(\Phi,\xi)$ without changing $V(\Phi)$ to make it holonomic near $V$.

\mn{\sc Step 5. Holonomic extension to $W$}.
Now on  $W\setminus \Op V$ we are in a position to apply 
Wrinkling Theorem 1.6B from \cite {[EM97]} (see also \cite {[EM98]}, p.335)
to extend the constructed $\phi_{A\cup V}$ as a leafwise wrinkled map
$\phi:(W,\F)\to\bbR$. The wrinkles of $\phi$ of any index have the 
canonical framing and thus this completes  the proof  of Theorem \ref{thm:main}.


\end{document}